
\input amstex
\loadbold
\documentstyle{amsppt}
\magnification=\magstep1
\TagsAsMath
\NoBlackBoxes



\pagewidth{6.5truein}
\pageheight{9.375truein}
\vcorrection{-0.375truein}


\long\def\ignore#1\endignore{#1}

\ignore 
\input xy \xyoption{matrix} \xyoption{arrow} 
\endignore


\def\im{\operatorname{Im}}
\def\rinj{\operatorname{inj}}
\def\id{\operatorname{id}}
\def\co{\operatorname{co}}
\def\mult{\underline{\operatorname m}}

\def\mmn{M_{m,n}}
\def\mmt{M_{m,t}}
\def\mtn{M_{t,n}}
\def\glt{GL_t}
\def\slt{SL_t}
\def\mmnlet{M_{m,n}^{\le t}}
\def\mmnleto{M_{m,n}^{\le t,\circ}}
\def\mnleto{M_n^{\le t,\circ}}
\def\mmto{\mmt^\circ}
\def\mtno{\mtn^\circ}

\def\O{{\Cal O}}
\def\ommn{\O(\mmn)}
\def\ommt{\O(\mmt)}
\def\omtn{\O(\mtn)}
\def\oglt{\O(\glt)}
\def\ommnlet{\O(\mmnlet)}
\def\ommnleto{\O(\mmnleto)}
\def\ommto{\O(\mmto)}
\def\omtno{\O(\mtno)}
\def\ov{\O(V)}
\def\ovo{\O(V^{\circ})}

\def\Oq{{\Cal O}_q}
\def\oqmn{\Oq(M_n)}
\def\oqmmn{\Oq(\mmn)}
\def\oqmmt{\Oq(\mmt)}
\def\oqmtn{\Oq(\mtn)}

\def\oqglt{\Oq(\glt)}
\def\oqslt{\Oq(\slt)}
\def\oqmmnlet{\Oq(\mmnlet)}

\def\oqmtno{\Oq(\mtno)}
\def\oqv{\Oq(V)}
\def\oqvmtn{\Oq(V_{m,t,n})}
\def\oqvo{\Oq(V^{\circ})}

\def\oqkst{\Oq(K^\times)}

\def\oqmnt{\Oq(M_{n,t})}
\def\oqmnto{\Oq(M_{n,t}^\circ)}
\def\oqmnlet{\Oq(M_n^{\le t})}
\def\oqmnleto{\Oq(M_n^{\le t,\circ})}

\def\oqgltchibar{\oqglt[\overline{\chi}]}

\def\roqst{\rho^*_q}
\def\lamqst{\lambda^*_q}
\def\muqst{\mu^*_q}
\def\gamqst{\gamma^*_q}
\def\Gamqst{\Gamma^*_q}
\def\thetqst{\theta^*_q}
\def\thqst{\theta_q^*}
\def\xiqst{\xi^*_q}

\def\iqst{i^*_q}
\def\jqst{j^*_q}

\def\phiqst{\phi_q^*}
\def\varphiqst{\varphi_q^*}

\def\roqsto{\rho^{*\circ}_q}
\def\lamqsto{\lambda^{*\circ}_q}
\def\muqsto{\mu^{*\circ}_q}
\def\gamqsto{\gamma^{*\circ}_q}

\def\A{{\Cal A}}

\def\C{{\Cal C}}
\def\I{{\Cal I}}

\def\NN{{\Bbb N}}
\def\CC{{\Bbb C}}
\def\ZZ{{\Bbb Z}}

\def\itil{\widetilde{I}}
\def\jtil{\widetilde{J}}

\def\robar{\overline{\rho}}
\def\rbar{\overline{r}}
\def\cbar{\overline{c}}
\def\chibar{\overline{\chi}}
\def\pbar{\overline{p}}


\def\DeCP{{\bf 1}}
\def\GoLen{{\bf 2}}
\def\Jor{{\bf 3}}
\def\LeSt{{\bf 4}}
\def\Mon{{\bf 5}}
\def\NYM{{\bf 6}}
\def\PaWa{{\bf 7}}
\def\Pro{{\bf 8}}

\topmatter

\title The First Fundamental Theorem of Coinvariant Theory
for the Quantum General Linear Group 
\endtitle

\rightheadtext{COINVARIANT THEORY FOR QUANTUM $GL_t$}

\author K. R. Goodearl, T. H. Lenagan, and L. Rigal \endauthor

\address Department of Mathematics, University of California, Santa
Barbara, CA 93106, USA\endaddress
\email goodearl\@math.ucsb.edu\endemail

\address Department of Mathematics, J.C.M.B., Kings Buildings,
Mayfield Road, Edinburgh EH9 3JZ, Scotland\endaddress
\email tom\@maths.ed.ac.uk\endemail

\address Universit\'e Jean Monnet
(Saint-\'Etienne), Facult\'e des Sciences et Techniques,
D\'e\-par\-te\-ment de Math\'ematiques, 23 rue du Docteur Paul
Michelon, 42023 Saint-\'Etienne C\'edex 2, France
\endaddress

\email Laurent.Rigal\@univ-st-etienne.fr\endemail

\abstract We prove First Fundamental Theorems of Coinvariant
Theory for the standard coactions of the quantum groups
$\Oq(GL_t(K))$ and
$\Oq(SL_t(K))$  on the quantized algebra $\Oq(M_{m,t}(K))
\otimes
\Oq(M_{t,n}(K))$.  (Here $K$ is an arbitrary field and $q$ an
arbitrary nonzero scalar.) In both cases, the set of
coinvariants is a subalgebra of $\Oq(M_{m,t}(K)) \otimes
\Oq(M_{t,n}(K))$, which we identify.
\endabstract

\subjclass 16W30, 17B37 \endsubjclass

\thanks The research of the first author was partially supported
by National Science Foundation research grant DMS-9622876, and
that of the first two authors by NATO Collaborative Research Grant
960250. The first author also thanks Mme\. M.-P\. Malliavin for
inviting him to the Universit\'e de Paris VI during January 1999,
where a portion of this work was done. Another portion was done
while the third author was attending the L.M.S\. Symposium on
Quantum Groups (Durham, July 1999) and during a subsequent visit to
the University of Edinburgh. He wishes to thank the London
Mathematical Society for partial financial support of both trips.
\endthanks

\endtopmatter

\document

\head Introduction\endhead

One of the highlights of classical invariant theory is the
determination of the algebra of invariant functions for the
standard action of the general linear group on the variety of
pairs of matrices over a field $K$. More precisely, the standard
action of $\glt= \glt(K)$ on the variety $V := \mmt(K)\times
\mtn(K)$ induces an action of $\glt$ on $\ov$, and the classical
theorem determines the algebra of invariants, $\ov^{\glt}$. We
recall the details below, since, if we assume that $K$ is
algebraically closed, the method of proof we follow has an easy
geometric translation.

The main theorem of this paper, Theorem 4.5, gives a quantum analog
of the above theorem. Since the quantum group $\oqglt$ is not a
group but a Hopf algebra, we first place the classical situation
into a Hopf algebra context. This is standard: the action of $\glt$
on $V$ induces a coaction of
$\oglt$ on $\ov$, under which $\ov$ becomes an $\oglt$-comodule
algebra, and $\ov^{\glt}$ equals the algebra of
$\oglt$-coinvariants in $\ov$. It is this situation which has a
natural quantization: the coordinate ring $\ov$ becomes the
algebra $\oqv := \oqmmt
\otimes \oqmtn$, and the coaction of $\oglt$ on $\ov$ becomes a
coaction of $\oqglt$ on $\oqv$. We prove the First Fundamental
Theorem of Coinvariant Theory for this coaction, that is, we
identify the set of coinvariants, $\oqv^{\co \oqglt}$.

There is a natural comultiplication map $\thqst:\oqmmn
\rightarrow
\oqmmt \otimes \oqmtn$, which is the quantum
analog of matrix multiplication $\mmt \times \mtn \rightarrow
\mmn$. We prove
that the set of coinvariants is equal to the image of $\thqst$. 
In an earlier paper, \cite{\GoLen}, the first two authors have
shown that the kernel of $\thqst$ is the ideal generated by the
$(t+1)\times (t+1)$  quantum minors of $\oqmmn$: this is the
Second Fundamental Theorem of Coinvariant Theory for this
comodule action.  Taken together, these two results give a
complete description of $\oqv^{\co\oqglt}$. Further, we
investigate the coaction of $\oqslt$ on $\oqv$ induced by that
of $\oqglt$ and identify the coinvariants of this coaction.

The basic structure of our proof follows the outline of one of
the possible proofs in the classical invariant theoretical
setting.   However, there are significant problems that arise
due to the noncommutative setting.  The most striking one is that,
unlike in the commutative case, the $\oqglt$-comodule $\oqv$ is not
a comodule algebra. For this reason, it is not even obvious at the
outset that the set of coinvariants forms a subalgebra. More
generally, the quantum analogs of several maps that we need are not
algebra morphisms, and so their properties cannot be analyzed
simply by checking how they behave on sets of algebra generators.

Nevertheless, it is useful to start by reviewing the classical
situation, to provide a skeleton for our approach.

\definition{The classical situation} 
 We fix an algebraically closed field $K$ and positive integers
$m,n,t$.  For integers $u,v >0$, we write
$M_{u,v}= M_{u,v}(K)$ for the set of $u \times v$ matrices with
entries in $K$. We will be mainly interested in the general
linear group
$\glt= \glt(K)$ and its standard action on the algebraic variety
$V= V_{m,t,n} :=\mmt \times \mtn$.   This action is given by:
\ignore
$$\xymatrixrowsep{0.25pc}\xymatrixcolsep{1.5pc}
\xymatrix{
\glt \times V \ar[r]^-{\gamma} & V \\
(g,(A,B)) \ar@{|->}[r] &(Ag^{-1},gB).
}$$ 
\endignore  
\noindent Thus $\glt$ acts on $\ov \cong \ommt \otimes
\omtn$. Classical invariant theory is interested in computing
the subalgebra $\ov^{\glt}$ of invariants for this action. The
description of this algebra goes as follows. Consider
the morphism of varieties
\ignore
$$\xymatrixrowsep{0.25pc}\xymatrixcolsep{2.0pc}
\xymatrix{
\mmt \times \mtn \ar[r]^-{\theta} & \mmn \\
(A,B) \ar@{|->}[r] &AB
}$$ 
\endignore  
\noindent and its associated comorphism
$\theta^*: \ommn \rightarrow \ommt \otimes \omtn$. Let $X_{ij}$
(for $1\le i\le m$ and $1\le j\le n$) stand for the usual
coordinate functions on the variety $\mmn$, and let
$\I_{t+1}$ denote the ideal of
$\ommn$ generated by all the $(t+1)\times (t+1)$ minors of the
generic matrix $(X_{ij})$ over $\ommn$. (This ideal is zero if
$t\ge \min\{m,n\}$.)
\enddefinition

\proclaim{Theorem 0.1} The ring of invariants $\ov^{\glt}$
equals $\im
\theta^*$. \qed
\endproclaim

\proclaim{Theorem 0.2} The kernel of $\theta^*$ is $\I_{t+1}$.
\qed
\endproclaim

Theorems 0.1 and 0.2 are respectively known as the {\it First
Fundamental Theorem of Invariant Theory\/} and the {\it Second
Fundamental Theorem of Invariant Theory\/} (for $\glt$). They
give a complete description of $\ov^{\glt}$.

We denote by
$\mmnlet$ the subvariety of
$\mmn$ of $m \times n$ matrices with rank at most $t$. This
variety is just the image of the morphism $\theta$, and so we can
factor $\theta$ in the form
$$\mmt\times\mtn @>{\mu}>> \mmnlet @>{\subseteq}>> \mmn.$$
Since the restriction map $r : \ommn\rightarrow\ommnlet$ is
surjective, the comorphism $\mu^*$ is injective and has the same
image as $\theta^*$. Thus, Theorem 0.1 can be rephrased in the
form $\ov^{\glt}= \im\mu^*$. Further, Theorem 0.2 shows that
$\ker(r)= \ker(\theta^*)= \I_{t+1}$, and therefore $\ommnlet=
\ommn/\I_{t+1}$.

\definition{The proof of De Concini and Procesi for Theorem 0.1}
We briefly describe the proof of Theorem 0.1 given by De Concini
and Procesi in \cite{\DeCP}; more precisely,  we follow the
exposition of that proof given in \cite{\Pro}.

The general case can be easily reduced to that where $t<
\min\{m,n\}$. So, we restrict attention to that particular case.
Let us fix some notation. We will denote by
$\mmnleto$ the open subset of $\mmnlet$ consisting of matrices with
rank at most
$t$ whose upper leftmost $t \times t$ minor is invertible. In a
similar way, $\mmto$ denotes the set of $m \times t$ matrices
whose uppermost $t
\times t$ minor is invertible, and $\mtno$ denotes the set of $t
\times n$ matrices whose leftmost $t
\times t$ minor is invertible. Finally, we set $V^{\circ} = \mmto
\times \mtno$. Clearly the action $\gamma$ of $\glt$ on
$V$ restricts to an action $\gamma^{\circ}$ on $V^{\circ}$ and we
can ask for the invariants of that restricted action. It turns
out that they are easy to compute.

Let us adopt the convention that if $(A,B)$ is an element of $V =
\mmt \times \mtn$ we write $A$ and $B$ in the forms
$$A= \left[ \matrix A_0 \\ A_1 \endmatrix \right]
\qquad\text{and}\qquad B= \left[ \matrix B_0 & B_1
\endmatrix\right]$$
 with $A_0,B_0 \in M_{t}$ while $A_1 \in
M_{m-t,t}$ and
$B_1 \in M_{t,n-t}$. 

   With these notations, we define a morphism of varieties:
\ignore
$$\xymatrixrowsep{0.25pc}\xymatrixcolsep{2.0pc}
\xymatrix{
V^{\circ} \ar[r]^-{i} & \mmnleto \times
\glt
\\ (A,B) \ar@{|->}[r] &(AB,B_0).
}$$ 
\endignore
\noindent It turns out that $i$ is actually an isomorphism of
varieties (this is an elementary fact, and the inverse morphism
can be explicitly written down) and that the  action of $\glt$ on
$\mmnleto \times \glt$ induced by
$\gamma^{\circ}$ is just the natural one:
\ignore
$$\xymatrixrowsep{0.25pc}\xymatrixcolsep{2.0pc}
\xymatrix{
\glt \times (\mmnleto \times \glt) \ar[r]^-{\xi} &
\mmnleto \times \glt \\ 
(g,C,h) \ar@{|->}[r] &(C,gh).
}$$ 
\endignore 
\noindent The invariants for $\xi$ are easy to compute,
using for instance the coinvariants of the associated
$\oglt$-comodule structure on $\O(\mmnleto \times \glt)$ (see
the quantum case below). One then proves that $\O(\mmnleto
\times \glt)^{\glt}=\ommnleto
\subseteq \O(\mmnleto \times \glt)$.

   To recover from this the invariants for $\gamma^{\circ}$ it is
enough to use the comorphism $i^{\ast}$ of $i$. Let us consider
the morphism 
$$\mu^{\circ} : \mmto \times \mtno = V^{\circ}
\longrightarrow \mmnleto$$
given by multiplication of matrices,  and denote by
$(\mu^{\circ})^{\ast}$ its comorphism:
$$ (\mu^{\circ})^{\ast} : \ommnleto
\longrightarrow \ovo=\ommto
\otimes \omtno.$$ 
 Then we can describe $i^{\ast}$ as the
composition
$$\align i^{\ast} :\ \ommnleto \otimes \oglt 
 & @>{(\mu^{\circ})^{\ast}\otimes \id}>>  
\ovo \otimes \oglt @>{\id \otimes {\subseteq}}>>
 \ovo \otimes \omtno \\ 
 & @>=>> \ommto \otimes \omtno \otimes \omtno @>{\id \otimes
\mult}>> \ovo, \endalign$$ 
where $\mult$ denotes multiplication
in the algebra $\omtno$.  From the above we get that the ring of
invariants for the action
$\gamma^{\circ}$ is   just $i^{\ast}(\ommnleto\otimes 1)=
(\mu^{\circ})^{\ast}(\ommnleto)$. 

The ring of invariants for
$\gamma^{\circ}$ can thus be described as the localisation of
$\im \mu^{\ast}$ with respect to the multiplicative set
generated by
$d_{Y} \otimes d_{Z}$ where $d_{Y}$ is
the uppermost $t
\times t$ minor of the generic matrix
$(Y_{ij})$ of generators of
$\ommt$, and $d_{Z}$ is
the leftmost $t
\times t$ minor of the generic matrix $(Z_{ij})$ of generators of
$\omtn$.

The second step of the proof is to show that one can ``remove
denominators'' to deduce invariants for the action $\gamma$ from
invariants for the localised action
$\gamma^\circ$. Indeed, if
$\phi$ is an invariant function in
$\ov$ (i.e., for the action
$\gamma$) then of course it is an invariant function in the
localised ring $\ovo=\ov_{d_{Y} \otimes d_{Z}}$ (for the action
$\gamma^{\circ}$). Hence, there is a non-negative integer $s$
such that $\phi (d_{Y} \otimes d_{Z})^{s}
\in
\im \mu^{\ast}$. Thus, in order to establish that $\ov^{\glt} =
\im \mu^{\ast}$, it is enough to prove that for $\psi \in \ov$,
if there is a non-negative integer $s$ such that $\psi (d_{Y}
\otimes d_{Z})^{s} \in
\im \mu^{\ast}$, then $\psi \in \im \mu^{\ast}$. This last
result is proved using the theory of standard bases.
\enddefinition

In the quantum situation, the lower left quantum minors play a
special r\^{o}le, since they are normal elements.  It is for
this reason that, unlike in the commutative case,   we will
invert leftmost {\it lower} 
$t\times t$ minors instead of leftmost upper ones.
\medskip 

Throughout the paper, we work over an arbitrary base field $K$
and make an arbitrary choice of a nonzero element $q\in K$.  We
will have to deal with the following four quantized coordinate
rings:
$\oqmmn$, $\oqmmt$, $\oqmtn$  and $\oqglt$. In order to avoid
confusion, we will denote their respective canonical generators
by
$X_{ij}$, $Y_{ij}$, $Z_{ij}$ and
$T_{ij}$.  The
definitions of these algebras will be as in \cite{\NYM, \GoLen};
for instance, $X_{ij}X_{ik}= qX_{ik}X_{ij}$ when $j<k$. Thus, each
time we use results from \cite{\PaWa} we must replace
$q$ by
$q^{-1}$.
 Finally, a convention concerning notation: each time
that we have to deal with the multiplication map in an algebra,
we will denote it by $\mult$. The context will make clear which
algebra is concerned.

\head 1. The general setup\endhead

Below, we will have to deal with the following situation combining a
right and a left comodule algebra. Let
$(H,\mult,\eta,\Delta,\varepsilon,S)$ be a Hopf algebra,
$(A,\rho^\ast)$ a right comodule algebra over $H$, and
$(B,\lambda^\ast)$ a left comodule algebra over $H$. Here, 
$$\rho^{\ast} : A \longrightarrow A \otimes H
\qquad\text{and}\qquad
\lambda^\ast : B \longrightarrow H \otimes B$$ 
are the comodule
structure maps; that $A$ and $B$ are {\it comodule algebras\/}
means that $\rho^*$ and $\lambda^*$ are also algebra morphisms. It
is well known that
$A$ can be turned into a left
$H$-co\-mod\-ule using the structure map 
$A @>{\rho^{\ast}}>> A \otimes H 
@>{\id \otimes S}>> A \otimes H
@>{\tau_{12}}>> H \otimes A$, where $\tau_{12}$ is the flip. Thus, $A
\otimes B$ can be equipped with the structure of a left
$H$-comodule  via the following structure map:
$$\align \gamma^\ast :  A \otimes B 
 & @>{\rho^{\ast} \otimes \lambda^\ast}>>   A
\otimes H \otimes H \otimes B
@>{\id \otimes S \otimes \id \otimes \id}>>  
A  
\otimes H \otimes H \otimes B \\
 & @>{\tau_{(132)}}>> H \otimes H
\otimes A \otimes B
@>{\mult \otimes \id \otimes \id}>>   H
\otimes A \otimes B, \endalign$$
where $\tau_{(132)}$ is the isomorphism which permutes the factors
according to the cycle $(132)$, that is, $\tau_{(132)}(a\otimes
h\otimes h'\otimes b)= h\otimes h'\otimes a\otimes b$. Using the
standard comodule notations
$\rho^\ast(a)=
\sum_{(a)} a_{0} \otimes a_{1}$ and  $\lambda^\ast(b)=
\sum_{(b)} b_{-1} \otimes b_{0}$ for $a \in A$ and
$b
\in B$ (cf.~\cite{\Mon, p.~11}), one thus has
$$
\gamma^\ast(a \otimes b) = 
\sum_{(a),(b)} S(a_{1})b_{-1} \otimes a_{0} \otimes b_{0}.$$ 
Of
course, $(A \otimes B,\gamma^\ast)$ is not a comodule algebra
any longer. Nevertheless, this comodule continues to have nice
properties that are now described. Let us recall
that if
$(M,\nu^\ast)$ is a left $H$-comodule, then the set of {\it
$H$-coinvariants of $(M,\nu^\ast)$\/} (or {\it
$\nu^\ast$-co\-in\-var\-i\-ants\/} for short) is the sub-comodule
of
$M$ defined by
$M^{\co H}:=\{x \in M \mid \nu^\ast(x)=1 \otimes x\}$. It is
immediate that if
$(M,\nu^\ast)$ is a comodule algebra then  $M^{\co H}$ is a
subalgebra of $M$. In the more general situation described above
this property is not automatic but still true.

\proclaim{Proposition 1.1} In the above notation: 

{\rm (a)} If $v \in A \otimes B$ is such that $\gamma^*(v)=
z\otimes v$ for some central element $z\in H$, and if 
$w \in A \otimes B$ is any element, then $\gamma^\ast(vw)
=\gamma^\ast(v)\gamma^\ast(w)$. In particular, this holds when $v$ is
a $\gamma^\ast$-coinvariant.

{\rm (b)} If $v,w \in A \otimes B$ are
$\gamma^\ast$-coinvariants, then
$vw$ is again a $\gamma^\ast$-coinvariant.

{\rm (c)} The set $(A \otimes
B)^{\co H}$ is a subalgebra of $A \otimes B$.
\endproclaim

\demo{Proof} Without loss of generality, we may assume that
$w=a' \otimes b'$ is a pure tensor. Moreover, let us write $v =
\sum_{i=1}^{r} a_{i} \otimes b_{i}$. Since both
$\rho^\ast$ and $\lambda^\ast$ are algebra morphisms, we have
$$(\rho^\ast \otimes \lambda^\ast)(vw) = \sum_{i=1}^{r}
\sum_{(a_{i}),(b_{i})} \sum_{(a'),(b')} a_{i,0}a^\prime_{0}
\otimes a_{i,1}a^\prime_{1} \otimes b_{i,-1}b^\prime_{-1}
\otimes b_{i,0} b^\prime_{0},
$$ 
 where $a_{i,0}= (a_i)_0$ etc. Hence,
$$\multline (\gamma^\ast)(vw) = \sum_{i=1}^{r}
\sum_{(a_{i}),(b_{i})}  \sum_{(a'),(b')}
S(a_{i,1}a^\prime_{1}) b_{i,-1}b^\prime_{-1}
\otimes a_{i,0}a^\prime_{0} 
\otimes b_{i,0} b^\prime_{0} \\
\aligned &=
\sum_{i=1}^{r}
\sum_{(a_{i}),(b_{i})}  \sum_{(a'),(b')} 
S(a^\prime_{1})S(a_{i,1}) b_{i,-1}b^\prime_{-1}
\otimes a_{i,0}a^\prime_{0} 
\otimes b_{i,0} b^\prime_{0} \\
 &=
\sum_{(a'),(b')}   \bigl(S(a^\prime_{1})\otimes 1 \otimes 1\bigr)
\biggl(\sum_{i=1}^{r}
\sum_{(a_{i}),(b_{i})} S(a_{i,1}) b_{i,-1}
\otimes a_{i,0} 
\otimes b_{i,0}\biggr) \bigl(b^\prime_{-1} \otimes a^\prime_{0}
\otimes b^\prime_{0}\bigr) \\
 &=
\sum_{(a'),(b')}   \bigl(S(a^\prime_{1})\otimes 1 \otimes 1\bigr)
\gamma^*(v) \bigl(b^\prime_{-1} \otimes a^\prime_{0} \otimes
b^\prime_{0}\bigr). 
\endaligned \endmultline$$ 
But, since $\gamma^*(v)= z\otimes v$ and $z$ is central, we
have 
$$\align  (\gamma^\ast)(vw)  &= \sum_{(a'),(b')}  
\bigl(S(a^\prime_{1})\otimes 1 \otimes 1\bigr) (z \otimes v)
\bigl(b^\prime_{-1}
\otimes a^\prime_{0} \otimes b^\prime_{0}\bigr) \\ 
&= (z \otimes v)
\sum_{(a'),(b')}  S(a^\prime_{1})b^\prime_{-1} \otimes a^\prime_{0}
\otimes b^\prime_{0} = \gamma^*(v) \gamma^\ast(w). 
\endalign$$ 
This proves statement (a). Clearly, (b) as well as
(c) follow at once from (a). \qed\enddemo

We now record a special case of a coaction for which computing the
coinvariants is easy.  Again,
$(H,\mult,\eta,\Delta,\varepsilon,S)$ is a Hopf algebra, and $M$
denotes a vector space. Clearly, the map 
\ignore
$$\xymatrixrowsep{0.25pc}\xymatrixcolsep{1.5pc}
\xymatrix{
M \ar[r] &H\otimes M\\
a \ar@{|->}[r] &1\otimes a
}$$
\endignore
\noindent defines a (trivial) left
coaction of
$H$ on $M$. Moreover, the comultiplication
$\Delta$ of $H$ makes $H$ itself into a left $H$-comodule. The
tensor product of these two coactions gives a left coaction of $H$
on $M \otimes H$:
\ignore
$$\xymatrixrowsep{0.25pc}\xymatrixcolsep{1.5pc}
\xymatrix{
M \otimes H \ar[r]^-{\xi^*} & H \otimes M \otimes H\\
a \otimes h \ar@{|->}[r] & \sum_{(h)} h_1 \otimes a \otimes h_2.
}$$
\endignore

\proclaim{Lemma 1.2} With the above
notation,
$(M \otimes H)^{\co H}=M \otimes 1$.
\endproclaim

\demo{Proof} Consider $x \in (M \otimes H)^{\co H}$, and write
$x=\sum_{i=1}^{s} a_{i} \otimes h_{i}$ where the $a_{i}$ are
linearly independent. Now
$$\sum_{i=1}^{s} \sum_{(h_i)} h_{i,1} \otimes a_i \otimes
h_{i,2} =
\xi^*(x) = 1\otimes x = \sum_{i=1}^{s} 1 \otimes a_i \otimes h_i.$$
Since the $a_i$ are linearly independent, we obtain
$$\Delta(h_i)= \sum_{(h_i)} h_{i,1} \otimes h_{i,2} =
1\otimes h_i$$ 
for each $i$, whence $h_i= (\id \otimes \varepsilon)
\Delta(h_i) = 1{\cdot}\varepsilon(h_i)$. This shows that the $h_i$
are scalars, and so $x\in M\otimes 1$. Thus, $(M \otimes H)^{\co H}
\subseteq M \otimes 1$. The reverse inclusion is obvious.
\qed\enddemo

\head 2. Quantization of the standard action of $\glt$ on $V$
\endhead

As in the introduction, we fix our base field $K$ and positive
integers $m,n,t$. By analogy with the commutative situation, we
put
$$  \oqv= \oqvmtn := \oqmmt \otimes \oqmtn.
$$ 
It is easily checked that $\oqv$ is an iterated skew
polynomial extension of the base field $K$, and so $\oqv$ is a
noetherian domain.  A quantum analogue of the action
$\gamma$ is obtained as a coaction of $\oqglt$ on $\oqv$ that
we now describe.

It is easy to check that one can define morphisms of algebras 
satisfying the following rules:
\ignore
$$\xymatrixrowsep{0.075pc}\xymatrixcolsep{2.0pc}
\xymatrix{
\oqmmt  \ar[r]^-{\roqst} &  \oqmmt \otimes \oqglt &\oqmtn
\ar[r]^-{\lamqst} & \oqglt \otimes \oqmtn \\
Y_{ij} \ar@{|->}[r] & \sum_{k=1}^{t} Y_{ik} \otimes T_{kj}
&Z_{ij} \ar@{|->}[r] & \sum_{k=1}^{t} T_{ik} \otimes
Z_{kj}.    
}$$   
\endignore
\noindent Moreover, $\roqst$ endows 
$\oqmmt$ with a right $\oqglt$-comodule
algebra structure, while $\lamqst$ endows $\oqmtn$ with a left
$\oqglt$-comodule algebra structure. It follows that
$\tau_{12} \circ (\id \otimes S)
\circ \roqst$ gives $\oqmmt$ the structure
of a left $\oqglt$-comodule, where $S$ denotes the antipode in
$\oqglt$. The tensor product of these two left coactions allows us to
define a left comodule structure on
$\oqv$ that we denote by $\gamqst$. So, as in Section 1,
$$
\gamqst \, : \, \oqv \longrightarrow \oqglt \otimes \oqv
$$ 
is given by the rule
$$
\gamqst(a \otimes b) =
\sum_{(a),(b)}  S(a_{1})b_{-1} \otimes a_{0} \otimes b_{0}
$$
for
$a \in \oqmmt$ and $b \in \oqmtn$.

The main objective of this paper is to calculate $\oqv^{\co \oqglt}$,
the set of coinvariants of $\gamqst$. As we see by Proposition
1.1, this set is a subalgebra of
$\oqv$.

\definition{Quantum minors} Recall that quantum minors in a
quantum matrix algebra $\Oq(M_{u,v})$ correspond to quantum
determinants in subalgebras of $\Oq(M_{u,v})$. More precisely, if
$1\le r_1< r_2< \dots< r_l\le u$ and $1\le c_1< c_2< \dots< c_l\le
v$, then the subalgebra of $\Oq(M_{u,v})$ generated by $\{
X_{r_ic_j} \mid 1\le i,j\le l\}$ is naturally isomorphic to
$\Oq(M_l)$, and the element of this subalgebra corresponding to
the quantum determinant in $\Oq(M_l)$ is called the quantum minor
of $\Oq(M_{u,v})$ corresponding to the rows $r_1,\dots,r_l$ and
columns $c_1,\dots,c_l$ (see \cite{\NYM, \S1.2; \PaWa, \S4.3}).
Since we require a number of the formulas and results developed
in
\cite{\GoLen}, we shall use the notation of that paper for quantum
minors. Thus, the $l\times l$ quantum minor described above will
be denoted by the symbol
$$[r_l r_{l-1} \cdots r_1 \mid c_1c_2 \cdots c_l],$$
or by $[R|C]$ where $R$ and $C$ denote the sets
$\{r_1,\dots,r_l\}$ and $\{c_1,\dots,c_l\}$, respectively.
\enddefinition

\definition{Remark 2.1} (i) The ideal of $\oqmmn$ generated by all
the $(t+1) \times (t+1)$ quantum minors will be denoted
$\I_{t+1}$, or
$\I^{m,n}_{t+1}$ when $m$ and $n$ require emphasis. (In case
$t\ge
\min\{m,n\}$, there are no $(t+1) \times
(t+1)$ quantum minors in $\oqmmn$, and $\I_{t+1}=0$.) Moreover, we
put
$$\oqmmnlet := \oqmmn/\I^{m,n}_{t+1}$$
and $x_{ij}:=X_{ij} + \I_{t+1} \in \oqmmnlet$.

(ii) For any $m'\ge m$ and $n'\ge n$, we will identify $\oqmmn$
with the subalgebra of $\Oq(M_{m',n'})$ generated by those $X_{ij}$
with $i\le m$ and $j\le n$. Now there is a
$K$-algebra retraction
$$\pi_{m,n}= \pi_{m,n}^{m',n'} : \Oq(M_{m',n'}) \longrightarrow
\oqmmn$$
such that $\pi_{m,n}(X_{ij})= X_{ij}$ for $1\le i\le m$ and $1\le
j\le n$, and $\pi_{m,n}(X_{ij})=0$ for $i>m$ or $j>n$. It is worth
mentioning here that the inclusion map $\oqmmn \rightarrow
\Oq(M_{m',n'})$ (corresponding to the identification above)
provides a section for $\pi_{m,n}$; this fact will be used without
further comment. Note that
$\pi_{m,n}[I|J]= [I|J]$ when $I\subseteq \{1,\dots,m\}$ and
$J\subseteq \{1,\dots,n\}$, while $\pi_{m,n}[I|J]=0$ otherwise. It
follows that $\pi_{m,n}(\I^{m',n'}_{t+1})= \I^{m,n}_{t+1}$, and
consequently $\I^{m,n}_{t+1}= \I^{m',n'}_{t+1} \cap \oqmmn$.

(iii) As is easily checked, there is a morphism
of algebras
\ignore
$$\xymatrixrowsep{0.25pc}\xymatrixcolsep{1.5pc}
\xymatrix{
\oqmmn  \ar[r]^-{\thetqst} &  \oqv \\
X_{ij} \ar@{|->}[r] & \sum_{k=1}^{t} Y_{ik} \otimes Z_{kj}.
}$$
\endignore
\noindent For instance, this can be checked by using the
commutative diagram

\ignore
$$\xymatrixrowsep{2.0pc}\xymatrixcolsep{4.0pc}
\xymatrix{
\Oq(M_l) \ar[r]^-{\Delta} &\Oq(M_l) \otimes \Oq(M_l) \ar[d]^{\pi_{m,t} \otimes
\pi_{t,n}}\\
\oqmmn \ar[r]^-{\thetqst} \ar[u]^{\subseteq} &\oqv
}$$
\endignore

\noindent where $l\ge \max\{m,n\}$ and $\Delta$ denotes the
comultiplication of $\Oq(M_l)$.

(iv) The comultiplication rule for quantum
minors in a square quantum matrix algebra, say $\Oq(M_l)$,
states that
$$\Delta [I|J]= \sum \Sb L\subseteq \{1,\dots,l\}\\ |L|=|I| \endSb
[I|L] \otimes [L|J]$$
for all quantum minors $[I|J]\in \Oq(M_l)$ (e.g., \cite{\NYM,
(1.9)}).

(v) From points (ii), (iii) and (iv), it is clear that $\I_{t+1}$ is
contained in the kernel of
$\thetqst$. Hence, there is an induced morphism of algebras
$$\muqst : \oqmmnlet
\longrightarrow \oqv$$  
such that $\muqst(x_{ij})=\sum_{k=1}^{t} Y_{ik} \otimes
Z_{kj}$. The work of the first
and second authors in \cite{\GoLen} leads to the following result:
\enddefinition

\proclaim{Theorem 2.2} In the above notation, the morphism
$\muqst$ is an injection. In particular, $\oqmmnlet$ is a domain.
\endproclaim

\demo{Proof} We must show that $\ker(\thetqst)= \I_{t+1}$. In view
of Remark 2.1(ii), there is no loss of generality in assuming
that $m=n>t$. That case is proved in \cite{\GoLen, Proposition
2.4}. \qed\enddemo

Clearly, Theorem 2.2 gives a quantum analogue of the Second
Fundamental Theorem of Invariant Theory; see Theorem 0.2.
The rest of this work is devoted to proving a quantum analogue of
the First Fundamental Theorem of Invariant Theory; see Theorem
0.1. A first easy proposition is in order.

\proclaim{Proposition 2.3} The set of coinvariants $\oqv^{\co
\oqglt}$ is a subalgebra of
$\oqv$ containing $\im
\muqst$.
\endproclaim
 
\demo{Proof} We are in the setup defined at the begining of
Section 1, so we may apply Proposition
1.1(c) to see that 
$\oqv^{\co \oqglt}$ is a subalgebra of
$\oqv$. For the rest of the statement, since the map
$\muqst$ is a morphism of algebras, it is enough to prove
that
$\muqst(x_{ij})$ is a
$\gamqst$-coinvariant, for $1 \le i \le m$ and $1 \le j \le n$.
However,
$$\align
\gamqst(\muqst(x_{ij}))  &=
\sum_{k=1}^{t}\gamqst(Y_{ik} \otimes Z_{kj})  =
\sum_{r,s,k=1}^{t} S(T_{rk})T_{ks} \otimes Y_{ir} \otimes
Z_{sj}\\ 
 &=
\sum_{r,s=1}^{t} \varepsilon(T_{rs}) \otimes Y_{ir} \otimes Z_{sj} =
\sum_{r=1}^{t} 1 \otimes Y_{ir} \otimes Z_{rj} = 1 \otimes
\muqst(x_{ij}). \qquad\square \endalign
$$ 
\enddemo

\definition{Remark 2.4} By using the retraction maps
$\pi_{\bullet,\bullet}$ discussed in Remark 2.1(ii), most results
proved for quantized coordinate rings of square matrices
immediately carry over to the rectangular case, as in the proof of
Theorem 2.2. For the purposes of the present paper, rather than
carry over a large number of results in this way, it is more
efficient to check that it suffices to prove our main theorem in
the case $m=n$. This also allows us to assume that $t<n$, as
follows.

Choose an integer $l\ge \max\{m,n\}$ (later, it will be helpful to
take $l> \max\{m,n\}$), and observe that we have a commutative
diagram as follows:

\ignore
$$\xymatrixrowsep{2.4pc}\xymatrixcolsep{8.0pc}
\xymatrix@!0{
\oqmmnlet \ar[dd]_{\muqst} \ar[dr]_{\subseteq} \ar@{=}[rr]
&&\oqmmnlet \ar[dd]^{\muqst}\\
 &\Oq(M_l^{\le t}) \ar[dd]^(0.3){\muqst}
\ar[ur]_{\hphantom{m}\overline{\pi}_{m,n}}\\
\oqvmtn \ar[dd]_{\gamqst} \ar[dr]_{\subseteq} \ar@{=}[rr]|\hole
&&\oqvmtn \ar[dd]^{\gamqst}\\
 &\Oq(V_{l,t,l}) \ar[dd]^(0.3){\gamqst}
\ar[ur]_{\hphantom{m}\pi_{m,t}\otimes\pi_{t,n}}\\
\oqglt\otimes\oqvmtn \ar[dr]_(0.4){\id\otimes\subseteq\hphantom{m}}
\ar@{=}[rr]|\hole &&\oqglt\otimes\oqvmtn\\
 &\oqglt\otimes\Oq(V_{l,t,l})
\ar[ur]_(0.6){\hphantom{mddd}\id\otimes\pi_{m,t}\otimes\pi_{t,n}}
}$$
\endignore

\noindent Here `$\subseteq$' is shorthand for various embeddings
induced by identifications (cf\. Remark 2.1(ii)), and
$\overline{\pi}_{m,n}$ denotes the map induced by the retraction
$\pi_{m,n}: \Oq(M_l)
\rightarrow
\oqmmn$ (recall from Remark 2.1(ii) that $\pi_{m,n}(\I^{l,l}_{t+1})=
\I^{m,n}_{t+1}$). It is clear from the above diagram that
$$\Oq(V_{l,t,l})^{\co\oqglt}= \muqst\bigl( \Oq(M_l^{\le
t}) \bigr) \quad\implies\quad \Oq(V_{m,t,n})^{\co\oqglt}=
\muqst\bigl( \oqmmnlet \bigr).$$
Hence, to prove our main theorem it suffices to consider the case
$(m=n>t)$.
\enddefinition

\head 3. A localisation of $\gamqst$\endhead

In view of Remark 2.4, we assume, until further notice, that
$m=n>t$. Certain notations contract slightly in this case; e.g.,
$\oqmmnlet$ becomes $\oqmnlet$.

\definition{Remark 3.1} (i)
We label the following $t\times t$ quantum minors:
$$\align d_{X} &:= [n,\dots,n-t+1 \mid 1,\dots,t] \in \oqmn\\
d_{Y} &:= [n,\dots,n-t+1 \mid 1,\dots,t] \in \oqmnt\\
d_{Z} &:= [t,\dots,1 \mid 1,\dots,t] \in \oqmtn\\
d_T &:= [t,\dots,1 \mid 1,\dots,t] \in \oqglt. \endalign$$
Further, we set $d_x= d_X+ \I_{t+1} \in \oqmnlet$. It is known
that $d_X\notin \I_{t+1}$ (for instance, it is clear from Remark
3.2(i) below that $\thetqst(d_X)\ne 0$), and so $d_x\ne 0$.

(ii) It is well known that in any $\Oq(M_{u,v})$, the ``lower
left'' quantum minors -- i.e., those corresponding to the lowest
$l$ rows and the leftmost $l$ columns, for any $l$ -- are normal
elements (e.g., see \cite{\GoLen, Corollary 5.2}). In particular,
the quantum minors $d_X$, $d_Y$, $d_Z$ are normal in $\oqmn$,
$\oqmnt$,
$\oqmtn$, respectively. Since $d_T$ is the quantum determinant in
$\oqglt$, it is central. It follows from the normality of $d_X$ that
$d_x$ is normal in $\oqmnlet$.

(iii) Since $d_Y$ is a normal
element in $\oqmnt$, we
can consider the localisation of $\oqmnt$ with respect to
the multiplicative set generated by $d_{Y}$. By analogy with the
commutative case, we denote this localisation,
$\oqmnt[d_Y^{-1}]$, by $\oqmnto$. Similarly, the localisations
$\oqmtn[d_Z^{-1}]$ and $\oqmnlet[d_{x}^{-1}]$ will be denoted
$\oqmtno$ and 
$\oqmnleto$, respectively. At this stage it is worth mentioning
that 
$$\oqmnleto\cong \bigl( \oqmn[d_X^{-1}] \bigr)
{\bigm/} \bigl( \I_{t+1}[d_X^{-1}] \bigr);$$
 this is an easy
consequence of the universal properties of quotients and
localisations. We will use this isomorphism without further
comment.
\enddefinition

\definition{Remark 3.2} (i) It follows from Remark 2.1(ii,iv)
that in the notation of Section
2, one has
$$\roqst(d_{Y})=d_{Y} \otimes d_{T} \qquad\text{ and }\qquad
\lamqst(d_{Z})=d_{T} \otimes d_{Z}.$$  
Further, $\thetqst(d_X)= d_Y\otimes d_Z$, and so
$$\muqst(d_{x})=d_{Y} \otimes d_{Z}.$$

(ii) In view of Remark 3.1(iii) and point (i) above, $\roqst$ and
$\lamqst$ extend uniquely to morphisms of algebras
$$\align  \roqsto &: \oqmnto
\longrightarrow  \oqmnto
\otimes \oqglt\\
 \lamqsto &: \oqmtno \longrightarrow
\oqglt \otimes \oqmtno. \endalign $$   

(iii) Clearly, $\lamqsto$ defines a left $\oqglt$-comodule
algebra structure on
$\oqmtno$, while
$\roqsto$ defines a right $\oqglt$-comodule algebra structure on
$\oqmnto$. Again, the composition $\tau_{12} \circ (\id \otimes S)
\circ  \roqsto$ defines a left
$\oqglt$-comodule structure on
$\oqmnto$. Finally, if we set
$$\oqvo:= \oqmnto \otimes \oqmtno$$
and tensor the above two left coactions, we obtain a map 
$$  \gamqsto \, : \, \oqvo
\longrightarrow \oqglt
\otimes \oqvo,
$$  
endowing $\oqvo$ with the structure of a left
$\oqglt$-comodule. The map $\gamqsto$ follows the same formula as
$\gamqst$; namely,
$$\gamqsto(a \otimes b) =
\sum_{(a),(b)}  S(a_{1})b_{-1} \otimes a_{0} \otimes b_{0}$$
for $a\in \oqmnto$ and $b\in \oqmtno$.
Obviously, the restriction of
$\gamqsto$ to $\oqv$ is just
$\gamqst$.

(iv) Recall from Theorem 2.2 that the morphism
$\muqst : \oqmnlet
\longrightarrow \oqv$ is injective. Because of Remark 3.1(iii) and
point (i) above, we can extend
$\muqst$ to the corresponding localised algebras and obtain
another injective morphism of algebras:
$$  \muqsto : \oqmnleto
\longrightarrow \oqvo.
$$  
\enddefinition

   We shall need a technical lemma. We have already mentioned
that
$\gamqst$ is not an algebra morphism, and so neither is $\gamqsto$.
Nevertheless, we have the following statement, the first part
of which is a special case of Proposition 1.1.

\proclaim{Lemma 3.3} The set of coinvariants $\oqvo^{\co \oqglt}$
is a subalgebra of $\oqvo$, containing $\im \muqsto$. 
\endproclaim

\demo{Proof} As noted, Proposition 1.1(c) implies that $\oqvo^{\co
\oqglt}$ is a subalgebra of $\oqvo$.

 Let us observe that for $r \in \NN$, the
element
$(d_{Y}
\otimes d_{Z})^{-r}$ is a coinvariant for
$\gamqsto$. Indeed, since both 
$\roqsto$ and $\lamqsto$ are
morphism of algebras,  we have 
$\roqsto (d_{Y}^{-r}) = d_{Y}^{-r} \otimes
d_{T}^{-r}$ and
$\lamqsto (d_{Z}^{-r}) = d_{T}^{-r} \otimes 
d_{Z}^{-r}$. Recall that $S(d_T)= d_T^{-1}$ \cite{\NYM, p\.
40; \PaWa, proof of Theorem 5.3.2}. Hence, 
$$\align \gamqsto((d_{Y} \otimes d_{Z})^{-r}) &=
\gamqsto(d_{Y}^{-r} \otimes d_{Z}^{-r}) = S(d_{T}^{-r})d_{T}^{-r}
\otimes  d_{Y}^{-r} \otimes d_{Z}^{-r}\\
 &= d_{T}^{r}d_{T}^{-r}
\otimes  d_{Y}^{-r} \otimes d_{Z}^{-r}  =  1
\otimes  d_{Y}^{-r} \otimes d_{Z}^{-r}.  \endalign$$ 
This proves the claim. Now since $\muqsto$ is a morphism of
algebras, any element $c\in \im \muqsto$ has the form
$c= a(d_{Y} \otimes d_{Z})^{-r}$ for some $a \in \im\muqst$ and
$r
\in \NN$. So, from the observation above, the fact that $\im
\muqst
\subseteq \oqv^{\co  \oqglt}$ (see Proposition 2.3), and the fact
that 
$\oqvo^{\co
\oqglt}$ is a subalgebra of $\oqvo$, we conclude that $c\in \oqvo^{\co
\oqglt}$. \qed\enddemo

\definition{Coaction of $\oqglt$ on $\oqmnleto \otimes \oqglt$}
One key point in the commutative case is that the action of $\glt$
on $V^{\circ}$  turns out to identify with the natural action of
$\glt$ on $\mnleto \times \glt$ by left
translation on the second factor. Here, we define and study a
quantum analogue for the latter action. 

We can apply the situation described at the end of Section
1 with  $M= \oqmnleto$ and
$H=\oqglt$. We thus obtain a left
$\oqglt$-comodule structure on 
$\oqmnleto \otimes \oqglt$ via:
\ignore
$$\xymatrixrowsep{0.25pc}\xymatrixcolsep{2.5pc}
\xymatrix{
\oqmnleto \otimes \oqglt \ar[r]^-{\xiqst} & \oqglt \otimes \oqmnleto \otimes
\oqglt \\  
x \otimes h \ar@{|->}[r] & \sum_{(h)} h_1 \otimes x \otimes
h_2.
}$$  
\endignore 
\noindent According to what we proved in Lemma 1.2, we have the
following statement:
\enddefinition

\proclaim{Proposition 3.4} The set of
$\xiqst$-coinvariants is a subalgebra of $\oqmnleto \otimes
\oqglt$. More precisely,
$$  \bigl( \oqmnleto \otimes \oqglt \bigr)^{\co \oqglt} = \oqmnleto
\otimes 1. \qquad\square
$$
\endproclaim

\definition{Link between the coactions $\gamqsto$
and
$\xiqst$} In light of the commutative case,
one expects that
$\oqmnleto \otimes \oqglt$
and $\oqvo$ are isomorphic as $\oqglt$-comodules. This is, in fact,
true, although we do not need the full statement here, only that
$\oqvo$ is a comodule retraction of $\oqmnleto \otimes \oqglt$. The
morphism that appears as the most natural candidate for an
isomorphism between these two comodules is a quantum analog of $i^*$
that we describe below.
\enddefinition

\definition{Remark 3.5}  (i) The
following injective homomorphisms of algebras will be of constant
use in the sequel. One should pay attention to the fact that, due to
the choice of the minors we have inverted, these are not the
usual maps one might expect to use. These embeddings are all into
the lower left corner of the target algebra, in order to make use
of the normality of the elements $d_Y$ etc. We denote all these
morphisms by ``${\rinj}$'' in order to avoid introducing too much
notation; the context will make clear which one we are using.

\ignore
$$\xymatrixrowsep{0.4pc}\xymatrixcolsep{1.5pc}
\xymatrix{
 &T_{ij} \ar@{|->}[rr] &&Y_{n-t+i,j}\\
T_{ij} \ar@{|->}[dd] &\oqglt \ar[rr]^{\rinj} \ar[dd]_{\rinj} &&\oqmnto
\ar[dd]^{\rinj}  &Y_{ij} \ar@{|->}[dd]\\ \\
Z_{ij} &\oqmtno \ar[rr]^{\rinj} &&\oqmnleto &x_{ij}\\
 &Z_{ij} \ar@{|->}[rr] &&x_{n-t+i,j}
}$$
\endignore
 
\noindent The injectivity of these maps is easy to prove, since
in each case we can exhibit a left inverse. We show how to deal
with the case of the left hand one; similar arguments hold for the
others.  The retraction
$\pi_{t,t}$ (see Remark 2.1(ii)) induces a morphism of
algebras                        
$\oqmtno \longrightarrow \oqglt$
that maps 
$Z_{ij}$ to $T_{ij}$ if $j \le t$ and to $0$ if $j >t$.  Clearly,
this map provides a left inverse for
$\rinj$.

(ii) Identify $\oqmnto$ (respectively,
$\oqmtno$) with its image in 
$\oqmnleto$. Then, the restriction of 
$\muqsto$ to $\oqmnto$ (respectively,  $\oqmtno$) coincides with
$\roqsto$ (respectively,
$\lamqsto$), as one can easily check on generators (this is
enough since the maps involved are morphism of algebras). More
precisely, we are claiming that the following diagram is
commutative:

\ignore
$$\xymatrixrowsep{1.75pc}\xymatrixcolsep{2.4pc}
\xymatrix{
\oqmnto \ar[r]^-{\rinj} \ar[d]_{\roqsto} &\oqmnleto 
\ar[d]^{\muqsto} &\oqmtno
\ar[l]_-{\rinj} \ar[d]^{\lamqsto}\\
\oqmnto \otimes \oqglt \ar[r]^-{\id \otimes \rinj} &\oqmnto
\otimes \oqmtno &\oqglt \otimes \oqmtno \ar[l]_-{\rinj \otimes \id}
}$$
\endignore

\noindent For example, to see that the left hand square commutes,
it suffices to check that $\muqsto\circ\rinj$ and $(\id\otimes
\rinj)\circ \roqsto$ agree on the $Y_{ij}$ and on $d_Y$. This
is clear -- both maps send $Y_{ij}$ to $\sum_{k=1}^t Y_{ik}
\otimes Z_{kj}$ and $d_Y$ to $d_Y \otimes d_Z$. The
commutativity of the right hand square is proved analogously.

(iii) Finally, the restriction of 
$\lamqsto$ (respectively, $\muqsto$) to $\oqglt$ coincides with
$\Delta$, i.e., the following diagram is commutative:

\ignore
$$\xymatrixrowsep{1.75pc}\xymatrixcolsep{2.4pc}
\xymatrix{
\oqglt \ar[r]^-{\rinj} \ar[d]_{\Delta} &\oqmtno \ar[r]^-{\rinj}
\ar[d]^{\lamqsto} &\oqmnleto \ar[d]^{\muqsto}\\ 
\oqglt \otimes \oqglt \ar[r]^-{\id \otimes \rinj} &\oqglt
\otimes \oqmtno \ar[r]^-{\rinj \otimes \id} &\oqmnto \otimes 
\oqmtno 
}$$
\endignore

In the following, we will identify the algebras $\oqglt$,
$\oqmnto$, $\oqmtno$ with their images under $\rinj$ and we will
use the commutativity of the above diagrams without explicit
mention.
\enddefinition

\definition{Quantization of $i^*$}
Due to our identification of $\oqglt$ with a subalgebra of
$\oqmtno$, we can use the multiplication map $\mult$ on
$\oqmtno$ to endow
$\oqmtno$ with the structure of a right 
$\oqglt$-module. We now consider the following quantum analogue
of $i^{\ast}$:
$$\align \iqst : \oqmnleto \otimes \oqglt & @>{\muqsto \otimes
\id}>>
\oqmnto \otimes \oqmtno \otimes \oqglt \\
 & @>{\id \otimes \mult}>>  \oqmnto \otimes \oqmtno = \oqvo.
\endalign$$
It can be proved that $\iqst$ is an isomorphism of left
$\oqglt$-comodules, but we shall not need the full result -- for
our purposes, it suffices to construct a right inverse for
$\iqst$.

Recall that the antipode $S$ in $\oqglt$ is bijective (e.g.,
\cite{\PaWa, Theorem 5.4.2}).

Consider the map $j_q^\ast : \oqvo \rightarrow
\oqmnleto \otimes \oqglt$ defined as the composition of the
following maps:
$$\align \oqvo &= \oqmnto \otimes \oqmtno @>{\roqsto \otimes
\lamqsto}>> \oqmnto \otimes \oqglt \otimes \oqglt \otimes
\oqmtno\\
 & @>{\id\otimes\id\otimes\Delta\otimes\id}>> \oqmnto \otimes
\oqglt \otimes \oqglt \otimes \oqglt \otimes \oqmtno\\
 & @>{\id\otimes S\otimes\id\otimes S^{-1}\otimes\id}>> \oqmnto
\otimes \oqglt \otimes \oqglt \otimes \oqglt \otimes \oqmtno\\
 & @>{\tau_{(135)(24)}}>> \oqmtno \otimes \oqglt \otimes \oqmnto
\otimes \oqglt \otimes \oqglt\\
 & @>{\rinj\otimes \rinj\otimes \rinj\otimes \mult}>> \oqmnleto
\otimes \oqmnleto \otimes \oqmnleto \otimes \oqglt\\
 & @>{\mult(\mult\otimes\id) \otimes\id}>> \oqmnleto \otimes
\oqglt. \endalign$$
Here $\tau_{(135)(24)}$ denotes the isomorphism which shuffles
the factors via the permutation $(135)(24)$, that is,
$$\tau_{(135)(24)}( x_1\otimes x_2\otimes x_3\otimes x_4\otimes
x_5)= x_5\otimes x_4\otimes x_1\otimes x_2\otimes
x_3.$$
 For a pure tensor
$a
\otimes b
\in
\oqvo$, we can give the following formula for $\jqst(a\otimes
b)$ using the notation 
$\roqst(a)= \sum_{(a)} a_{0} \otimes a_{1}$ and $(\Delta \otimes
\id) \lamqst(b)  = \sum_{(b)} b_{-2} \otimes b_{-1}
\otimes b_{0}$. Namely,
$$j_{q}^\ast(a \otimes b) =
\sum_{(a),(b)} b_{0}S^{-1}(b_{-1})a_{0} \otimes S(a_{1})b_{-2}.
$$
\enddefinition

\proclaim{Lemma 3.6} The map $\jqst$ is a morphism of left
$\oqglt$-comodules. \endproclaim

\demo{Proof} We must prove that the following
diagram is commutative:

\ignore
$$\xymatrixrowsep{1.75pc}\xymatrixcolsep{3.0pc}
\xymatrix{
 \oqvo \ar[r]^-{\jqst} \ar[d]_{\gamqsto} &\oqmnleto \otimes \oqglt
\ar[d]^{\xiqst}\\
\oqglt \otimes \oqvo \ar[r]^-{\id \otimes
\jqst} &\oqglt \otimes \oqmnleto \otimes \oqglt
}$$
\endignore

\noindent If $a\in \oqmnto$ and $b\in \oqmtno$, then (with the
usual conventions for expanding components)
$$\align (\id\otimes\jqst) \gamqsto (a\otimes b) &=
\sum_{(a),(b)} S(a_1)b_{-1} \otimes \jqst(a_0\otimes b_0)\\
 &= \sum_{(a),(b)} S(a_2)b_{-3} \otimes b_0S^{-1}(b_{-1})a_0
\otimes S(a_1)b_{-2}. \endalign$$ 
Since $\Delta(S(c)d)= \sum_{(c),(d)} S(c_2)d_1
\otimes S(c_1)d_2$ for $c,d\in \oqglt$, we also have
$$\align \xiqst\jqst(a\otimes b) &= \sum_{(a),(b)} \xiqst \bigl(
b_0S^{-1}(b_{-1})a_0 \otimes S(a_1)b_{-2} \bigr)\\
 &= \sum_{(a),(b)}  S(a_2)b_{-3} \otimes b_0S^{-1}(b_{-1})a_0
\otimes S(a_1)b_{-2}. \endalign$$ 
Therefore $(\id\otimes\jqst)
\gamqsto (a\otimes b)=
\xiqst\jqst(a\otimes b)$, as desired. \qed\enddemo

In the sequel, we shall need the identity
$$\sum_{(c)} c_2S^{-1}(c_1)= \sum_{(c)} S^{-1}(c_2)c_1=
\varepsilon(c)\cdot 1$$
for $c\in \oqglt$. This follows from the defining property of
the antipode $S$ by applying the antiautomorphism $S^{-1}$.
Moreover, in view of the identity $\Delta\circ S=
\tau_{12}\circ (S\otimes S)\circ \Delta$ (e.g., \cite{\Mon,
Proposition 1.5.10}), we also have
$$(S^{-1} \otimes
S^{-1}) \circ
\tau_{12} \circ \Delta = \Delta \circ S^{-1}.$$

\proclaim{Lemma 3.7} The composition $\iqst\circ\jqst$ is
the identity map on $\oqvo$.
\endproclaim

\demo{Proof} Let $a\otimes b\in \oqvo$. Using Remark 3.5(ii,iii)
and the fact that
$\muqsto$ is a morphism of algebras, we have 
$$\align (\muqsto \otimes \id ) j_{q}^\ast(a \otimes b) &=
\sum_{(a),(b)} \bigl( \muqsto(b_{0}) \bigl(
\muqsto S^{-1}(b_{-1}) \bigr) \muqsto(a_{0}) \bigr) \otimes
S(a_{1})b_{-2}\\
 &= \sum_{(a),(b)} \bigl( \lamqsto(b_{0}) \bigl( \Delta
S^{-1}(b_{-1}) \bigr) \roqsto(a_{0}) \bigr) \otimes
S(a_{1})b_{-2}\\
 &= \sum_{(a),(b)} \bigl( (b_{-1}\otimes b_0) \bigl( \Delta
S^{-1}(b_{-2}) \bigr) (a_0\otimes a_1) \bigr) \otimes
S(a_{2})b_{-3}. \endalign$$ 
Consequently, with the help of the identities above, we see that
$$\align (\muqsto \otimes \id ) \jqst(a \otimes b) &=
\sum_{(a),(b)}  b_{-1} S^{-1}(b_{-2}) a_{0} 
\otimes b_{0} S^{-1}(b_{-3}) a_{1} 
\otimes S(a_{2})b_{-4} \\
 &= \sum_{(a),(b)} 
\varepsilon(b_{-1}) a_{0} 
    \otimes b_{0} S^{-1}(b_{-2}) a_{1} 
    \otimes S(a_{2})b_{-3} \\
 &= \sum_{(a),(b)}   a_{0} 
    \otimes b_{0} S^{-1} \bigl( b_{-2}\varepsilon(b_{-1}) \bigr)
a_{1} 
    \otimes S(a_{2})b_{-3} \\
 &= \sum_{(a),(b)}   a_{0} 
    \otimes b_{0} S^{-1}(b_{-1}) a_{1} 
    \otimes S(a_{2})b_{-2}. \endalign$$ 
Thus,
$$\align \iqst j_{q}^\ast(a \otimes b)  &= 
\sum_{(a),(b)}  a_{0} \otimes b_{0} S^{-1}(b_{-1}) a_{1}
S(a_{2})b_{-2} = \sum_{(a),(b)}   a_{0} \otimes b_{0}
S^{-1}(b_{-1}) \varepsilon(a_{1}) b_{-2} \\
 &= \sum_{(a),(b)}  
a_{0} \varepsilon(a_{1}) \otimes b_{0} S^{-1}(b_{-1}) b_{-2} = 
\sum_{(b)} a \otimes b_{0}\varepsilon(b_{-1}) = a \otimes b.
\qquad\square \endalign$$
\enddemo

\proclaim{Proposition 3.8} $\oqvo^{\co \oqglt} =\im \muqsto$.
\endproclaim

\demo{Proof} We have the inclusion ``$\supseteq$'' by Lemma 3.3.
To prove the reverse inclusion, consider an arbitrary
$\gamqsto$-coinvariant $c\in \oqvo$. Since $\jqst$ is a comodule
morphism (Lemma 3.6), $\jqst(c)$ is a $\xiqst$-coinvariant in
$\oqmnleto \otimes \oqglt$. Hence, $\jqst(c) \in \oqmnleto
\otimes 1$ by Proposition 3.4. Since $\iqst\jqst(c)=c$ by Lemma
3.7, and since clearly $\iqst(\oqmnleto \otimes 1)= \im
\muqsto$, we conclude that $c\in \im\muqsto$. This verifies the
inclusion ``$\subseteq$''. \qed\enddemo

\head 4. The First Fundamental Theorem for
$\oqglt$-coinvariants\endhead

In this section, we complete the proof of our main theorem --
that the set of $\gamqst$-co\-in\-var\-i\-ants in $\oqv$ equals
the image of
$\muqst$. We first discuss how to obtain $\gamqst$-coinvariants
from
$\gamqsto$-coinvariants by ``removing denominators''.

\definition{Remark 4.1} (i) As we already
mentioned in Remark 3.2(iii), the restriction of
$\gamqsto$ to
$\oqv$ is  just $\gamqst$. So, if $w \in \oqv$ is a
$\gamqst$-coinvariant, it is also a
$\gamqsto$-coinvariant. Thus, according to Proposition
3.8,
$w\in  \im \muqsto$. This means that there exists $s \in \NN$
such that
$w (d_{Y} \otimes d_{Z})^s \in \im \muqst$.
   So, to prove that  the set of $\gamqst$-coinvariants equals
$\im \muqst$, it is enough to show the following property
(cf.~Proposition 2.3):
$$\forall w \in \oqv, \; \text{ if } \; w (d_{Y} \otimes d_{Z})
\in \im \muqst,\; \text{ then }\; w \in \im \muqst. \tag\dagger
$$

(ii) Recall from  Remark 3.2(i) that the image under $\muqst$ of
the (normal) element $d_x \in \oqmnlet$ is the
(normal) element $d_{Y} \otimes d_{Z} \in \oqv$. We claim that
to prove $(\dagger)$, it is enough (in fact, equivalent) to
verify the following property:
$$(\muqst)^{-1} \bigl( \langle d_{Y}
\otimes d_{Z}\rangle \bigr)=  \langle d_x\rangle. \tag\ddagger$$
Thus, assume that $(\ddagger)$ holds, and let $w\in \oqv$ such
that $w (d_{Y} \otimes d_{Z})
\in \im \muqst$. Then $w (d_{Y} \otimes d_{Z})= \muqst(y)$ for
some $y\in \oqmnlet$, and $(\ddagger)$ implies that $y= zd_x$
for some $z\in \oqmnlet$. Since $\muqst$ is a morphism of
algebras, $w (d_{Y} \otimes d_{Z})= \muqst(y)= \muqst(z) (d_{Y}
\otimes d_{Z})$, and so $w= \muqst(z)$ because $\oqv$ is a
domain. This shows that $(\ddagger)$ implies $(\dagger)$. The
converse follows easily from the injectivity of $\muqst$ (see
Theorem 2.2).
\enddefinition

\definition{Remark 4.2} (i) The computations in this section and
in Section 6 below rely on the {\it preferred bases\/} for the
algebras
$\oqmn$,
$\oqmnt$, and
$\oqmtn$ developed in \cite{\GoLen, Section 1}, and we follow the
notation of that paper. See, in particular, \cite{\GoLen,
Corollary 1.11} for the rectangular case. We recall the notation
$[T|T']$ for the product of quantum minors corresponding to an
allowable bitableaux $(T,T')$, and
we recall also that it is sometimes convenient to label rows of
$(T,T')$ in the form $(I,J)$ where $I$ and $J$ are sets of row and
column indices, respectively.

(ii) On tensoring the preferred basis elements of $\oqmnt$ with
those of $\oqmtn$, we obtain a preferred basis for $\oqv$
consisting of all pure tensors $[S|S']\otimes [T|T']$ where
$(S,S')$ and $(T,T')$ are preferred bitableaux with appropriate
entries. In particular, $S$ and $T$ each have at most $t$ columns.

(iii) Because of the quantum Laplace expansions \cite{\NYM,
Proposition 1.1 and Corollary; \PaWa, Corollary 4.4.4}, the ideal
$\I_{t+1}$ in
$\oqmn$ contains all the
$r\times r$ quantum minors for $t+1\le r\le n$. Hence, $[T|T']\in
\I_{t+1}$ whenever $(T,T')$ is a preferred bitableau and $T$ has
at least $t+1$ columns.

The statement and proof of \cite{\GoLen, Proposition 2.4} show
that $\I_{t+1}$ does not contain any nonzero linear combinations
of distinct products $[S|S']$ where the $(S,S')$ are preferred
bitableaux and $S$ has at most $t$ columns. Consequently,
$\I_{t+1}$ has a $K$-basis consisting of the products $[T|T']$ such
that
$(T,T')$ is a preferred bitableau and $T$ has at least $t+1$
columns.

(iv) In view of point (iii), the cosets $[T|T'] +\I_{t+1}$ such
that  $(T,T')$
is a preferred bitableau and $T$ has at most $t$ columns form a
$K$-basis for $\oqmnlet$. To simplify the notation, we shall write
cosets $a+\I_{t+1}$ in the form $\overline a$. Thus,
$\overline{[T|T']} =0$ whenever $T$ has more than $t$ columns.
\enddefinition

It is convenient to label two special index sets:
$$\itil= \{n-t+1, \dots,n\} \qquad\quad\text{and}\qquad\quad
\jtil=
\{1,\dots,t\}.$$
With this notation, observe that $d_X$ and $d_Y$ can both be
labeled $[\itil|\jtil]$, while $d_Z$ and $d_T$ can both be
labeled $[\jtil|\jtil]$.

\proclaim{Lemma 4.3} {\rm (a)} A $K$-basis for $\langle
d_Y\rangle$ consists of all $[T|T']$ in $\oqmnt$ where $(T,T')$
is a preferred bitableau with first row $(\itil,\jtil)$.

{\rm (b)} A $K$-basis for $\langle
d_Z\rangle$ consists of all $[T|T']$ in $\oqmtn$ where $(T,T')$
is a preferred bitableau with first row $(\jtil,\jtil)$.

{\rm (c)} A $K$-basis for $\langle d_{Y}\otimes
d_{Z}\rangle$ consists of all pure tensors $[S|S'] \otimes
[T|T']$ in $\oqv$ where $(S,S')$ and $(T,T')$ are preferred
bitableaux with first rows $(\itil,\jtil)$ and $(\jtil,\jtil)$,
respectively.

{\rm (d)} A $K$-basis for $\langle
d_x\rangle$ consists of all $\overline{[T|T']}$ in $\oqmnlet$
where
$(T,T')$ is a preferred bitableau with first row $(\itil,\jtil)$.
\endproclaim

\demo{Proof} (a) Obviously any such $[T|T']$ lies in $\langle
d_Y\rangle$, because $d_Y= [\itil|\jtil]$. Recall that $d_Y$ is
a normal element of $\oqmnt$, so that $\langle d_Y\rangle =d_Y
\oqmnt$. Since $(\itil,\jtil)$ is the minimum index pair labeling
quantum minors in
$\oqmnt$, the preferred basis of
$\oqmnt$ is closed under left multiplication by
$d_Y$. Hence, $\langle d_Y\rangle$ is spanned by products
$d_Y[S|S']$ as $(S,S')$ runs through all preferred bitableaux,
and each such $d_Y[S|S']= [T|T']$ for some preferred bitableau
$(T,T')$ with first row
$(\itil,\jtil)$.

(b) This is proved in the same manner as part (a).

(c) Note that because $d_Y$ and $d_Z$ are normal
elements of $\oqmnt$ and $\oqmtn$, respectively,
$$\langle d_Y\otimes d_Z\rangle= (d_{Y}\otimes d_{Z}) \oqv=
d_Y\oqmnt \otimes d_Z\oqmtn= \langle d_Y\rangle \otimes \langle
d_Z\rangle.$$
Hence, part (c) follows directly from parts (a) and (b).

(d) By Remark 4.2(iv), a $K$-basis for $\oqmnlet$ consists of
all $\overline{[T|T']}$ such that $(T,T')$ is a preferred
bitableau and $T$ has at most $t$ columns. Among index pairs
$(I,J)$ with $|I|\le t$, the minimum element is $(\itil,\jtil)$.
Since
$d_x=
\overline{[\itil|\jtil]}$, we therefore obtain part (d) in exactly
the same manner as part (a). \qed\enddemo

In the following proof, we shall need the natural multigradings
on quantum matrix algebras (cf\. \cite{\GoLen, \S1.5}). For
instance,
$\oqmnt$ is graded by
$\ZZ^n
\times \ZZ^t$ with each generator $Y_{ij}$ having degree
$(\epsilon_i, \epsilon_j)$, where $\epsilon_1, \epsilon_2,
\dots$ denote the standard basis elements in $\ZZ^n$ and
$\ZZ^t$. In the proof, we use the label `homogeneous' to refer
to homogeneous elements with respect to the above gradings.

\proclaim{Proposition 4.4}  $(\muqst)^{-1} \bigl( \langle d_{Y}
\otimes d_{Z}\rangle \bigr)=  \langle d_x\rangle$. \endproclaim

\demo{Proof} The inclusion `$\supseteq$' is
clear since $\muqst(d_x)= d_Y\otimes d_Z$. If
this inclusion is proper, choose an element
$$x= \sum_{i=1}^r \alpha_i \overline{[T_i|T'_i]} \in
(\muqst)^{-1}\bigl( \langle d_Y\otimes d_Z\rangle \bigr)
\setminus \langle d_x\rangle$$
where the $\alpha_i$ are nonzero scalars, the $(T_i,T'_i)$ are
distinct preferred bitableaux, and the $T_i$ have at most $t$
columns. We may assume that none of the $\overline{[T_i|T'_i]}$ lie
in $\langle d_x\rangle$. Thus, by Lemma 4.3(d), none of the 
$(T_i,T'_i)$ has first row $(\itil,\jtil)$. Define
$n$-tuples
$\robar(T_i)$ as in
\cite{\GoLen, \S2.2}, and  let $\robar$ be the minimum of the
$\robar(T_i)$ under reverse lexicographic order. 

After re-indexing, we may assume that there is some $r'$ such that
$\robar(T_i)= \robar$ for $i\le r'$ and $\robar(T_i) >_{\text{rlex}}
\robar$ for $i>r'$. Applying \cite{\GoLen, Lemma 2.3} to each $\muqst
\overline{[T_i|T'_i]}$ ($= \thetqst [T_i|T'_i]$) and collecting
terms, we see (using the notation of
\cite{\GoLen, \S2.2}) that
$$\muqst(x)= \sum_{i=1}^{r'} \alpha_i [T_i|\mu(T_i)] \otimes
[\mu'(T_i)|T'_i] + \sum_j X_j\otimes Y_j$$
where the $X_j$ and $Y_j$ are homogeneous with $\cbar(X_j)=
\rbar(Y_j)
 >_{\text{rlex}} \robar$. We then observe (as in the proof of
\cite{\GoLen, Proposition 2.4}) that all of the $X_j$ belong to
different homogeneous components than the $[T_i|\mu(T_i)]$ for $i\le
r'$. Since $\muqst(x) \in \langle d_Y\otimes d_Z\rangle= \langle
d_Y\rangle \otimes \langle d_Z\rangle$ and the
ideal $\langle d_Y\rangle$ is homogeneous, it follows that
$$\sum_{i=1}^{r'} \alpha_i [T_i|\mu(T_i)] \otimes
[\mu'(T_i)|T'_i] \in \langle d_Y\otimes d_Z\rangle.$$

For $1\le i<j\le r'$, either $T_i\ne T_j$ or $T'_i\ne T'_j$,
whence
$(T_i,\mu(T_i)) \ne (T_j,\mu(T_j))$ or $(\mu'(T_i),T'_i) \ne
(\mu'(T_j),T'_j)$. Recall from \cite{\GoLen, \S2.2} that the
$(T_i,\mu(T_i))$ and the $(\mu'(T_i),T'_i)$ are preferred
bitableaux. In view of Lemma 4.3(c), it follows that
$T_i$ has first row $\itil$ and
$T'_i$ has first row $\jtil$ for
$1\le i\le r'$. However, this contradicts our choices above, and
therefore the proposition is proved. \qed\enddemo

\proclaim{Theorem 4.5} Let $m,n,t$ be arbitrary positive
integers. The set of
$\gamqst$-coinvariants in $\oqv= \oqmmt \otimes \oqmtn$ is $\im
\muqst$, that is, 
$$\bigl( \oqmmt \otimes \oqmtn \bigr)^{\co \oqglt}= \muqst\bigl(
\oqmmnlet \bigr).$$
\endproclaim

\demo{Proof} In view of Remark 2.4, there is no loss of generality
in assuming that $m=n>t$. The theorem is then immediate using 
Remark 4.1 and Proposition 4.4.
\qed\enddemo

\head 5. $\oqslt$-coinvariants\endhead

For arbitrary positive integers $m,n,t$, there is a natural coaction
$\Gamqst$  of
$\oqslt$ on
$\oqv$ induced from the
$\oqglt$-coaction $\gamqst$ that we have been studying. Denote by
$\pi$ the natural quotient map $\oqglt \rightarrow \oqslt$.  The
coaction
$\Gamqst: \oqv \rightarrow \oqslt \otimes \oqv$ is given by  
$\Gamqst:= (\pi \otimes \id)\circ \gamqst$.  In this section and
the next, we compute the $\Gamqst$-coinvariants. 

In fact, $\Gamqst$ can be constructed from the left coaction of
$\oqslt$ on $\oqmtn$ given by 
$$ (\pi \otimes \id)\circ \lamqst: \oqmtn \rightarrow \oqslt
\otimes \oqmtn
$$ and the right coaction of  $\oqslt$ on $\oqmmt$ given by  
$$ (\id \otimes \pi)\circ \roqst: \oqmmt \rightarrow \oqmmt
\otimes \oqslt
$$ by using the construction described in Section 1.  

\proclaim{Proposition 5.1} The set of coinvariants $\oqv^{\co
\oqslt}$ is a subalgebra of $\oqv$ containing $\im \muqst$.
\endproclaim

\demo{Proof} That $\oqv^{\co \oqslt}$ is a subalgebra
of $\oqv$ follows from Proposition 1.1(c). It contains
$\im \muqst$ because the elements of $\im \muqst$ are
$\gamqst$-coinvariants and hence also $\Gamqst$-coinvariants.
\qed\enddemo

In order to describe the $\Gamqst$-coinvariants, we may assume that
$m=n$ and that $t<n$, as in Remark 2.4.   We start by observing
that there is a natural $\ZZ$-grading on quantum matrices which
will simplify the problem.  The algebra $\oqmnt$ can be graded
by total degree in the variables $Y_{ij}$ and we put 
$
\oqmnt = \bigoplus_{i\ge 0} \oqmnt_i,
$  
where $\oqmnt_i$ is the subspace spanned by monomials of total
degree
$i$ in $\oqmnt$.  In the same way, and with obvious notation,
$
\oqmtn = \bigoplus_{j\ge 0} \oqmtn_j.
$ 
It follows that 
$$
\oqv = \bigoplus_{i,j=0}^\infty \oqv_{i,j},
$$ where $\oqv_{i,j} = \oqmnt_i \otimes \oqmtn_j$. 

Recall that $\roqst: \oqmnt \rightarrow \oqmnt \otimes \oqglt$ is
the algebra morphism defined by $\roqst(Y_{ij}) = \sum_{k=1}^t\,
Y_{ik}
\otimes T_{kj}$.  Thus, it is clear that 
$
\roqst(\oqmnt_i) \subseteq \oqmnt_i \otimes \oqglt.
$   
Similarly, one has 
$\lamqst(\oqmtn_j) \subseteq \oqglt \otimes \oqmtn_j.
$ It follows that 
$$ \align
\gamqst (\oqv_{i,j}) &\subseteq \oqglt \otimes \oqv_{i,j}\\
\Gamqst (\oqv_{i,j}) &\subseteq \oqslt \otimes \oqv_{i,j}
\endalign$$  
for all $i,j$. 

Given $w\in \oqv$, we will write $w= \sum_{i,j}\, w_{ij}$ with
$w_{ij}
\in \oqv_{i,j}$.  

\proclaim{Lemma 5.2}  Let $w= \sum_{i,j}\, w_{ij} \in \oqv$. 
Then 

{\rm (a)} $w$ is a $\gamqst$-coinvariant if and only if each
$w_{ij}
$ is a
$\gamqst$-coinvariant. 

{\rm (b)} $w$ is a $\Gamqst$-coinvariant if and only if each
$w_{ij}
$ is a 
$\Gamqst$-coinvariant.
\endproclaim

\demo{Proof}  If each $w_{ij} $ is a
$\gamqst$-coinvariant then, obviously, $w$ is a
$\gamqst$-coinvariant.  Conversely, suppose that $w$ is a
$\gamqst$-coinvariant. One has
$\sum_{i,j}\, 1\otimes w_{ij} = 1 \otimes w = \gamqst(w) =
\sum_{i,j}\,
\gamqst(w_{ij})$, and from the discussion above, we see that 
$\gamqst(w_{ij}) \in \oqglt \otimes \oqv_{i,j}$.  Thus, since 
$\oqglt \otimes \oqv = \bigoplus_{i,j}\, \oqglt \otimes
\oqv_{i,j}$, one has
$\gamqst(w_{ij}) = 1\otimes w_{ij}$ for all $i,j$.
This finishes the proof of (a); the proof of (b) is similar. 
\qed\enddemo 

The previous result
shows that, in order to describe the
$\Gamqst$-coinvariants, it is enough to describe the coinvariants
which are in each $\oqv_{i,j}$.  It is obvious from the definition
of $\Gamqst$ that any $\gamqst$-coinvariant will be a
$\Gamqst$-coinvariant.  Our first  aim  is to show that any
homogeneous
$\Gamqst$-coinvariant $w$ is a {\it
$\gamqst$-semi-coinvariant\/}, meaning that
$\gamqst(w) = d_T^s
\otimes w$ for some integer $s$; this is achieved in Theorem
5.5.   We begin by constructing some maps which will help us to
prove this fact.

The above grading on $\oqv$ relates to the following
comodule structure. Consider the Hopf algebra $\oqkst = K[T^{\pm
1}]$ and the algebra morphism
$
\pi': \oqglt \twoheadrightarrow \oqkst
$  such that $T_{ii} \mapsto T$ for all $i$, while $T_{ij}
\mapsto 0$ for
$i\neq j$. It is easily checked that $\pi'$ is, in fact, a
morphism of Hopf algebras. The composite map
$$
\oqv @>{\gamqst}>> \oqglt \otimes \oqv @>{\pi' \otimes \id}>>
\oqkst
\otimes \oqv$$
thus gives $\oqv$ the structure of a left comodule
over $\oqkst$.  The following proposition shows that, in fact,
$\oqv$ is an
$\oqkst$-comodule algebra. 

\proclaim{Proposition 5.3} {\rm (a)} $(\pi' \otimes
\id)\circ \gamqst $ is an algebra morphism.

{\rm (b)} $(\pi' \otimes \id) \gamqst(Y_{ij}\otimes 1)= 
T^{-1} \otimes Y_{ij} \otimes 1$ and $(\pi' \otimes \id)
\gamqst(1\otimes Z_{ij})= T \otimes 1\otimes Z_{ij}$ for all
$i,j$.

{\rm (c)} $(\pi' \otimes \id) \gamqst(w) = T^{j-i}\otimes w$ for
all $w\in \oqv_{i,j}$.
\endproclaim

\demo{Proof} (a) Recall that, for $a\in \oqmnt$ and $b\in
\oqmtn$, if we put
$\roqst(a) = \sum_{(a)}\,a_0 \otimes a_1$ and $\lamqst(b) =
\sum_{(b)}\, b_{-1} \otimes b_0$, then 
$
\gamqst(a\otimes b) =  \sum_{(a),(b)} S(a_1)b_{-1}\otimes a_0
\otimes b_0.$ 
	If, also, $a'\otimes b' \in \oqv$ then, since $\roqst$ and
$\lamqst$ are algebra morphisms, one has 
$$
\gamqst \bigl( (a\otimes b) (a'\otimes b') \bigr)  =
\sum_{(a'),(b')} S(a_1')S(a_1)b_{-1}b_{-1}' \otimes a_0a_0' \otimes
b_0b_0'.
$$  
Using the fact that $\oqkst$ is a commutative algebra we then
get 
$$ (\pi' \otimes \id) \gamqst \bigl( (a\otimes b) (a'\otimes b')
\bigr) = \bigl( (\pi'
\otimes \id) \gamqst \bigr) (a\otimes b) \bigl( (\pi' \otimes
\id) \gamqst \bigr) (a'\otimes b').
$$ 

(b) If $A_{jk}$ denotes the $(t-1)\times (t-1)$
quantum minor of $\oqglt$ obtained by deleting row $j$ and column
$k$, then $S(T_{kj}) = (-q)^{k-j} A_{jk}d_T^{-1}$, by \cite{\PaWa,
Theorem 5.3.2}, and so 
$$
\gamqst(Y_{ij} \otimes 1) = \sum_{k=1}^t \, S(T_{kj}) \otimes
Y_{ik}
\otimes 1 = \sum_{k=1}^t \, (-q)^{k-j}A_{jk}d_T^{-1} \otimes Y_{ik}
\otimes 1
$$ 
and 
$\gamqst(1\otimes Z_{ij}) = \sum_{k=1}^t\, T_{ik} \otimes 1 \otimes
Z_{kj}.
$  It follows that 
$(\pi' \otimes \id) \gamqst (Y_{ij} \otimes 1) = T^{-1}\otimes
Y_{ij} \otimes 1$  and that 
$(\pi' \otimes \id) \gamqst (1\otimes Z_{ij}) = T\otimes 1
\otimes Z_{ij}$.

(c) This is clear from parts (a) and (b). \qed\enddemo 

We shall need the algebra morphism 
$$ \phiqst: \oqglt \otimes \oqkst \rightarrow \oqgltchibar:=
\oqglt[\chi]/ \langle \chi^t - d_T^{-1}
\rangle,$$  
where $\chi$ is an indeterminate and $\chibar  = \chi + \langle
\chi^t - d_T^{-1} \rangle$, such that each
$T_{ij}\otimes 1 \mapsto T_{ij}\chibar $  and $1\otimes T
\mapsto
\chibar ^{-1}$.    It is easy to check that
$\phiqst(d_T\otimes 1) = d_T\chibar ^t =1$, and so there is
an induced algebra morphism 
$$\varphiqst: \oqslt \otimes \oqkst \rightarrow \oqgltchibar$$  
such that $\varphiqst \circ (\pi\otimes\id)= \phiqst$.

Finally, let $\alpha_q : \oqglt \rightarrow \oqgltchibar$ denote
the canonical injection.

\proclaim{Lemma 5.4} {\rm (a)} $\phiqst \circ(\id \otimes \pi')
\circ
\Delta = \alpha_q$.

{\rm (b)} $(\varphiqst \otimes \id)\circ (\id \otimes \pi'
\otimes \id) \circ (\id \otimes \gamqst) \circ \Gamqst =
(\alpha_q \otimes \id)\circ
\gamqst$.
\endproclaim 

\demo{Proof} (a) It is enough to check this on the generators
$T_{ij}$, since each of the maps involved is an algebra
morphism.  One has 
$$\phiqst (\id \otimes \pi') \Delta(T_{ij}) = \phiqst (\id
\otimes \pi') \bigl( \sum_{k=1}^t\, T_{ik}\otimes T_{kj} \bigr) =
\phiqst(T_{ij}
\otimes T) = T_{ij}.$$

(b) Since $\gamqst$ is the structure map of a left
$\oqglt$-comodule, we have $(\Delta\otimes\id) \circ\gamqst=
(\id\otimes\gamqst) \circ\gamqst$. In view of part (a), it
follows that
$$\align (\alpha_q \otimes \id)\circ \gamqst &=
(\phiqst\otimes\id) \circ (\id \otimes \pi' \otimes \id) \circ
(\Delta\otimes\id) \circ\gamqst\\
 &= (\phiqst\otimes\id) \circ (\id \otimes \pi'\otimes \id) \circ
(\id\otimes\gamqst) \circ\gamqst. \endalign$$
We also have the following commutative diagram:

\ignore
$$\xymatrixrowsep{3.2pc}\xymatrixcolsep{8.0pc}
\xymatrix@!0{
 &\oqv \ar[dl]_{\gamqst} \ar[dr]^{\Gamqst}&\\
\oqglt\otimes\oqv \ar[d]_{ \id\otimes \gamqst} \ar[rr]^{\pi\otimes
\id} & &\oqslt\otimes\oqv \ar[d]^{\id\otimes \gamqst}\\
\oqglt\otimes\oqglt\otimes\oqv \ar[d]_{\id\otimes\pi'\otimes \id}
\ar[rr]^{\pi\otimes\id\otimes \id} &&
\oqslt\otimes\oqglt\otimes\oqv \ar[d]^{\id \otimes \pi'\otimes \id}\\
\oqglt\otimes\oqkst\otimes\oqv \ar[rr]^{\pi\otimes \id\otimes \id}   
\ar[dr]_(0.4){\phi_q^* \otimes \id}&&
\oqslt\otimes\oqkst\otimes\oqv \ar[dl]^(0.4){\varphi_q^* \otimes \id}\\
 &\oqgltchibar\otimes \oqv
}$$
\endignore

\noindent Thus $(\phiqst\otimes\id) \circ (\id \otimes \pi'\otimes \id) \circ
(\id\otimes\gamqst) \circ\gamqst= (\varphiqst \otimes \id)\circ
(\id \otimes \pi'
\otimes \id) \circ (\id \otimes \gamqst) \circ \Gamqst$, which
completes the proof of part (b).
\qed\enddemo

\proclaim{Theorem 5.5}  Let $w\in \oqv_{i,j} $ be a
$\Gamqst$-coinvariant. 

{\rm (a)} The difference $i-j$ is divisible by $t$, and
$\gamqst(w) =d_T^{-s} \otimes w$ where $s= (i-j)/t$.

{\rm (b)} If $i\ge j$, then $(1\otimes d_Z^s)w$ is a
$\gamqst$-coinvariant.

{\rm (c)} If $i\le j$, then $(d_Y^{-s}\otimes 1)w$ is a
$\gamqst$-coinvariant. 
\endproclaim 

\demo{Proof} (a) By Proposition 5.3(c), we have
$(\pi'\otimes \id) \gamqst(w) = T^{j-i}\otimes w $.  From
this and the assumption that $\Gamqst(w)= 1\otimes w$, it follows
that 
$$\align (\varphiqst \otimes \id) (\id\otimes \pi'\otimes \id)
(\id\otimes \gamqst) \Gamqst(w) &= 
(\varphiqst \otimes \id) (\id\otimes \pi'\otimes \id)
(\id\otimes \gamqst) (1\otimes w)\\
&= (\varphiqst
\otimes \id)(1\otimes T^{j-i} \otimes w) = \chibar ^{i-j} \otimes
w. \endalign$$
Thus, by Lemma 5.4(b), $(\alpha_q\otimes\id) \gamqst(w)= \chibar
^{i-j} \otimes w$. In particular,
$$\chibar ^{i-j} \otimes w \in \im(\alpha_q\otimes\id)=
\oqglt\cdot 1 \otimes \oqv.$$
Since $\oqgltchibar$ is
a free
$\oglt$-module with basis $\{ 1, \chibar ,
\dots,
\chibar ^{t-1}\}$,  we have $t\mid i-j$. Set $s= (i-j)/t$.
Now
$$(\alpha_q\otimes\id) \gamqst(w)= \chibar
^{ts} \otimes w= (\alpha_q\otimes\id)(d_T^{-s} \otimes w),$$
and therefore $\gamqst (w) = d_T^{-s}
\otimes w$, as required.

(b) Since $\lamqst$ is an algebra morphism, we see using Remark
3.2(i) that $\lamqst(d_Z^s)= d_T^s\otimes d_Z^s$, and hence
$\gamqst(1\otimes d_Z^s)= d_T^s\otimes 1\otimes d_Z^s$. As
$d_T^s$ is central in $\oqglt$, we may apply Proposition 1.1(a)
to conclude that
$$\gamqst \bigl( (1\otimes d_Z^s)w \bigr)= \gamqst(1\otimes
d_Z^s) \gamqst(w)= (d_T^s\otimes 1\otimes d_Z^s)(d_T^{-s}
\otimes w)= 1\otimes \bigl( (1\otimes d_Z^s)w \bigr).$$

(c) By Remark 3.2(i), $\roqst(d_Y^{-s})= d_Y^{-s}\otimes
d_T^{-s}$. Since $S(d_T)= d_T^{-1}$, we thus obtain 
$$\gamqst(
d_Y^{-s}\otimes 1)= d_T^s\otimes d_Y^{-s}\otimes 1,$$
 and therefore
part (c) follows from a second application of Proposition
1.1(a).
 \qed\enddemo

\head 6. The First Fundamental Theorem for $\oqslt$-coinvariants
\endhead

In this final section, we determine the $\Gamqst$-coinvariants
in $\oqv$. Let $\A_1$ and $\A_2$ denote the respective
subalgebras of $\oqmmt$ and $\oqmtn$ generated by all the
$t\times t$ quantum minors. We shall prove that the set of
$\Gamqst$-coinvariants in $\oqv$ is the subalgebra generated by
$\A_1\otimes \A_2$ together with $\im \muqst$. In fact, this
subalgebra turns out to be equal to the product $(\A_1\otimes
\A_2)\cdot \im \muqst$. On the road to this goal, Theorem 5.5
puts us in a position roughly similar to that of Remark
4.1(ii), and to finish we need some computations analogous to
Proposition 4.4.

Until further notice, we continue to assume that $m=n>t$.

\definition{Remark 6.1} (i) As in Section 4, we set $\itil=
\{n-t+1,
\dots,n\}$ and
$\jtil= \{1,\dots,t\}$. Let
$P_1$ and $P_2$ denote the following ideals:
$$\align P_1 &= \big\langle\ \overline{[\itil|J]} \bigm| J\subseteq
\{1,\dots,n\}
\text{\ and\ } |J|=t\ \big\rangle \vartriangleleft \oqmnlet\\
P_2 &= \big\langle\ \overline{[I|\jtil]} \bigm| I\subseteq
\{1,\dots,n\}
\text{\ and\ } |I|=t\ \big\rangle \vartriangleleft \oqmnlet.
\endalign$$
These turn out to be completely prime ideals of $\oqmnlet$ --
see Remark 6.4. In view of \cite{\GoLen, Corollary 5.2}, the
quantum minors
$[\itil|J]$ can be arranged in a polynormal sequence, as can the
$[I|\jtil]$. Hence,
$$P_1= \sum \Sb J\subseteq \{1,\dots,n\}\\ |J|=t \endSb
\overline{[\itil|J]}
\oqmnlet \qquad\text{and}\qquad P_2= \sum \Sb I\subseteq
\{1,\dots,n\}\\ |I|=t \endSb \overline{[I|\jtil]} \oqmnlet.$$

(ii) It can be proved that a $K$-basis for $P_1$ consists of
all
$\overline{[T|T']}$ where $(T,T')$ is a preferred bitableaux with
first row of the form
$(\itil, ?)$, and that a $K$-basis for $P_2$ consists of all
$\overline{[T|T']}$ where $(T,T')$ is a preferred bitableaux with
first row of  the form
$(?,\jtil)$. However, we shall not need these bases.

(iii) In view of point (i), powers of the $P_i$ can be expressed
as follows:
$$\align P_1^s &= \sum \Sb J_1,\dots,J_s\subseteq \{1,\dots,n\}\\
|J_1|=\cdots= |J_s|=t \endSb \overline{ [\itil|J_1] [\itil|J_2]
\cdots [\itil|J_s]} \oqmnlet\\
P_2^s &= \sum \Sb I_1,\dots,I_s\subseteq \{1,\dots,n\}\\
|I_1|=\cdots= |I_s|=t \endSb \overline{ [I_1|\jtil] [I_2|\jtil]
\cdots [I_s|\jtil]} \oqmnlet \endalign$$
for $s=1,2,\dots$.
\enddefinition

\proclaim{Lemma 6.2} Let $s$ be a nonnegative integer.

{\rm (a)} A $K$-basis for $\langle d_Y^s\otimes 1\rangle$
consists of all pure tensors $[S|S']\otimes [T|T']$ in $\oqv$
where
$(S,S')$ and
$(T,T')$ are preferred bitableaux and the first $s$ rows of
$(S,S')$ equal $(\itil,\jtil)$.

{\rm (b)} A $K$-basis for $\langle 1\otimes d_Z^s \rangle$
consists of all pure tensors $[S|S']\otimes [T|T']$ in $\oqv$
where
$(S,S')$ and
$(T,T')$ are preferred bitableaux and the first $s$ rows of
$(T,T')$ equal $(\jtil,\jtil)$.
\endproclaim

\demo{Proof} Since $\langle d_Y^s\otimes 1\rangle= \langle
d_Y^s\rangle \otimes\oqmtn$ and $\langle 1\otimes d_Z^s \rangle=
\oqmnt\otimes \langle d_Z^s\rangle$, it suffices to check that
the ideals $\langle
d_Y^s\rangle$ and $\langle d_Z^s\rangle$ have bases of the
appropriate forms. This can be done just as in Lemma 4.3(a)(b).
\qed\enddemo

\proclaim{Proposition 6.3} If $s$ is any nonnegative integer,
then
$$(\muqst)^{-1} \bigl( \langle d_Y^s\otimes 1\rangle
\bigr)= P_1^s \qquad\quad \text{and} \qquad\quad (\muqst)^{-1}
\bigl(
\langle 1\otimes d_Z^s
\rangle
\bigr)= P_2^s.$$
\endproclaim

\demo{Proof} We prove the first equality; the proof of the second
is analogous. The proof closely mimics that of Proposition 4.4.

In view of Remark 2.1(iii-v), we see that
$$\muqst \overline{[\itil|J]}= [\itil|\jtil] \otimes [\jtil|J]
=d_Y \otimes [\jtil|J]$$
for all $t$-element subsets $J\subseteq \{1,\dots,n\}$. Using
Remark 6.1(iii) and the fact that $\muqst$ is an algebra
morphism, it follows that $P_1^s\subseteq (\muqst)^{-1} \bigl(
\langle d_Y^s\otimes 1\rangle
\bigr)$.

If
this inclusion is proper, choose an element
$$x= \sum_{i=1}^r \alpha_i \overline{[T_i|T'_i]} \in
(\muqst)^{-1}\bigl( \langle d_Y^s\otimes 1\rangle
\bigr) \setminus P_1^s$$
where the $\alpha_i$ are nonzero scalars, the $(T_i,T'_i)$ are
distinct preferred bitableaux, and the $T_i$ have at most $t$
columns. We may assume that none of the $\overline{[T_i|T'_i]}$ lie
in $P_1^s$. Let $\robar$ be the minimum of the
$\robar(T_i)$ under reverse lexicographic order. 

After re-indexing, we may assume that there is some $r'$ such that
$\robar(T_i)= \robar$ for $i\le r'$ and $\robar(T_i) >_{\text{rlex}}
\robar$ for $i>r'$. Applying \cite{\GoLen, Lemma 2.3} to each $\muqst
\overline{[T_i|T'_i]}$ and collecting terms, we see that
$$\muqst(x)= \sum_{i=1}^{r'} \alpha_i [T_i|\mu(T_i)] \otimes
[\mu'(T_i)|T'_i] + \sum_j X_j\otimes Y_j$$
where the $X_j$ and $Y_j$ are homogeneous with $\cbar(X_j)=
\rbar(Y_j)
 >_{\text{rlex}} \robar$. We then observe that all of the $X_j$
belong to different homogeneous components  than the
$[T_i|\mu(T_i)]$ for
$i\le r'$. Since 
$$\muqst(x) \in \langle d_Y^s\otimes 1\rangle=
\langle d_Y^s\rangle \otimes\oqmtn$$
and the ideal $\langle
d_Y^s\rangle$ is homogeneous, it follows that
$$\sum_{i=1}^{r'} \alpha_i [T_i|\mu(T_i)] \otimes
[\mu'(T_i)|T'_i] \in \langle d_Y^s\otimes 1\rangle.$$

For $1\le i<j\le r'$, either $(T_i,\mu(T_i)) \ne (T_j,\mu(T_j))$
or $(\mu'(T_i),T'_i) \ne (\mu'(T_j),T'_j)$. In view of Lemma
6.2(a), it follows that for
$1\le i\le r'$, the first $s$ rows of $T_i$ equal $\itil$. But
then $\overline{[T_i|T'_i]} \in P_1^s$ for $1\le i\le r'$, which
contradicts our choices above and therefore establishes the
desired equality.
\qed\enddemo

\definition{Remark 6.4} Although we shall not need the fact
here, we note that the case $s=1$ of Proposition 6.3 implies
that $P_1$ and $P_2$ are completely prime ideals of $\oqmnlet$.
To see this, it obviously suffices to show that $\langle
d_Y\otimes 1\rangle$ and $\langle 1\otimes d_Z
\rangle$ are completely prime ideals of $\oqv$. Observe that
$\oqmnt/\langle d_Y\rangle$ is an iterated skew polynomial ring
over an isomorphic copy of $\Oq(M_t)$ modulo the ideal generated
by its quantum determinant. Since the latter algebra is a domain
(\cite{\Jor}, \cite{\LeSt, p\. 182}, or see
\cite{\GoLen, Theorem 2.5}), it follows that
$\oqmnt/ \langle d_Y\rangle$ is a domain. Now the algebra
$$\oqv/\langle d_Y\otimes 1\rangle \cong \bigl( \oqmnt/
\langle d_Y\rangle \bigr) \otimes \oqmtn$$
is an iterated skew polynomial ring over $\oqmnt/ \langle
d_Y\rangle$, and therefore $\oqv/\langle
d_Y\otimes 1\rangle$ is a domain. Thus $\langle
d_Y\otimes 1\rangle$ is completely prime, and a similar argument
shows that $\langle 1\otimes d_Z
\rangle$ is completely prime.
\enddefinition

\proclaim{Lemma 6.5} Let $w\in \oqv_{i,j} $ be a
$\Gamqst$-coinvariant. 

{\rm (a)} If $i\ge j$, then $w\in (\A_1\otimes 1)\cdot
\im \muqst$.

{\rm (b)} If $i\le j$, then $w\in (1\otimes \A_2)\cdot
\im \muqst$.
\endproclaim

\demo{Proof} By Theorem 5.5(a), the ratio $s= (i-j)/t$ is an
integer.

(a) In this case, Theorem 5.5(b) shows that $(1\otimes d_Z^s)w$
is a $\gamqst$-coinvariant. Thus, by Theorem 4.5, $(1\otimes
d_Z^s)w= \muqst(x)$ for some $x\in \oqmnlet$. Now $\muqst(x)\in
\langle 1\otimes d_Z^s\rangle$, and hence $x\in P_2^s$, by
Proposition 6.3. In view of Remark 6.1(iii), we thus have $x=
\sum_{i=1}^r \pbar_ix_i$ for some $x_i\in \oqmnlet$ and some
$p_i\in \oqmn$ of the form
$$p_i= [I_{i,1}|\jtil] [I_{i,2}|\jtil] \cdots [I_{i,s}|\jtil].$$
Now by Remark 2.1(iii-v), each
$$\muqst \overline{[I_{i,l}|\jtil]}= [I_{i,l}|\jtil] \otimes
[\jtil|\jtil]= [I_{i,l}|\jtil] \otimes d_Z.$$
Since $\muqst$ is an algebra morphism, it follows that each
$$\muqst(\pbar_ix_i)= (1\otimes d_Z^s) \bigl( [I_{i,1}|\jtil]
\otimes 1 \bigr) \bigl( [I_{i,2}|\jtil] \otimes 1 \bigr)
\cdots\bigl( [I_{i,s}|\jtil] \otimes 1 \bigr)
\muqst(x_i).$$
Consequently, $(1\otimes d_Z^s)w= \muqst(x)= (1\otimes d_Z^s)v$
for some $v\in (\A_1\otimes 1)\cdot
\im \muqst$. Since $1\otimes d_Z$ is a nonzero element of the
domain $\oqv$, we conclude that $w=v$, and the proof of part (a)
is complete.

(b) This is proved in a similar fashion. \qed\enddemo

\proclaim{Theorem 6.6} Let $m,n,t$ be arbitrary positive
integers, and let $\A_1$ and $\A_2$ denote the respective
subalgebras of $\oqmmt$ and $\oqmtn$ generated by all the
$t\times t$ quantum minors. The set of $\Gamqst$-coinvariants in
$\oqv= \oqmmt\otimes\oqmtn$ is the subalgebra generated by
$\A_1\otimes \A_2$ and $\im \muqst$. More precisely,
$$\bigl( \oqmmt\otimes\oqmtn \bigr)^{\co \oqslt}= \bigl(
\A_1\otimes \A_2 \bigr) \cdot \muqst \bigl( \oqmmnlet \bigr).$$
\endproclaim

\demo{Proof} As in Remark 2.4, it is enough to consider the case
that $m=n>t$; we leave the details of the reduction to the
reader. Let $\C$ denote the set of $\Gamqst$-coinvariants in
$\oqv$, and recall from Proposition 5.1 that $\C$ is a
subalgebra of $\oqv$ containing $\im \muqst$.

If $[I|J]$ is a $t\times t$ quantum minor in
$\oqmnt$ (necessarily, $J=\jtil$), then as in Remark
2.1(iii)(iv), we see that
$\roqst [I|J]= [I|J] \otimes d_T$. It follows that $\gamqst(
[I|J]\otimes 1)= d_T^{-1} \otimes [I|J]\otimes 1$, whence
$$\Gamqst( [I|J]\otimes 1)= (\pi\otimes\id) (d_T^{-1} \otimes
[I|J]\otimes 1)= 1\otimes [I|J] \otimes 1$$
and so $[I|J]\otimes 1\in \C$. Similarly, $1\otimes [I|J] \in
\C$ for all $t\times t$ quantum minors $[I|J]$ in $\oqmtn$.
Therefore
$\C$ contains the subalgebra of $\oqv$ generated by $\A_1\otimes
\A_2$ and $\im \muqst$.

On the other hand, it is immediate from Lemmas 5.2(b) and
6.5 that $\C$ is contained in $\bigl(
\A_1\otimes \A_2 \bigr) \cdot \im \muqst$. The theorem follows.
\qed\enddemo

\enddocument